\documentclass[3p,times]{elsarticle}

\usepackage{ecrc}

\volume{00}

\firstpage{1}

\journalname{Procedia Computer Science}

\runauth{}

\jid{procs}

\jnltitlelogo{Procedia Computer Science}

\CopyrightLine{2011}{Published by Elsevier Ltd.}

\usepackage{amssymb}
\usepackage{amsthm}

\usepackage[figuresright]{rotating}

\usepackage{anyfontsize}
\usepackage[pagewise]{lineno}
\usepackage{epsfig,latexsym,amssymb, amsmath,amscd,amsxtra}
\usepackage{bbm}
\usepackage{subfigure}
\usepackage{indentfirst}
\usepackage{graphicx}
\usepackage{algorithm}  
\usepackage{algpseudocode}    
\usepackage{bbm}
\usepackage{multirow}
\usepackage{subfigure}
\usepackage{color}
\usepackage[numbers]{natbib}
\usepackage{indentfirst}
\usepackage{comment}
\usepackage[misc]{ifsym}
\newcommand{\RNum}[1]{\uppercase\expandafter{\romannumeral #1\relax}}
\newtheorem{remark}{Remark}

\begin{document}

\begin{frontmatter}




\title{Deep FBSDE Neural Networks for Solving Incompressible Navier-Stokes Equation and Cahn-Hilliard Equation}


\author[1]{Yangtao Deng}
\ead{ytdeng1998@foxmail.com}

\author[1]{Qiaolin He \corref{cor}}
\ead{qlhejenny@scu.edu.cn}

\cortext[cor]{Corresponding author}

\address[1]{School of Mathematics, Sichuan University, Chengdu, China}

\begin{abstract}
    Efficient algorithms for solving high-dimensional partial differential equations (PDEs) has been an exceedingly difficult task for a long time, due to the curse of dimensionality. We extend the forward-backward stochastic neural networks (FBSNNs) which depends on forward-backward stochastic differential equation (FBSDE)  to solve incompressible Navier-Stokes equation. For Cahn-Hilliard equation, we derive a modified Cahn-Hilliard equation from a widely used stabilized scheme for original Cahn-Hilliard equation. This equation can be written as  a continuous parabolic system, where FBSDE can be applied and the unknown solution is approximated by neural network. Also our method is successfully developed to Cahn-Hilliard-Navier-Stokes (CHNS) equation. The accuracy and stability of our methods are shown in many numerical experiments, specially in high dimension.
\end{abstract}


\end{frontmatter}



\section{Introduction}\label{sec01}
\setcounter{equation}{0}
High-dimensional nonlinear partial differential equations (PDEs) are used widely in a number of areas of social and natural sciences. Due to the significant nonlinearity of nonlinear PDEs, particularly in high-dimensional cases, analytical solutions to nonlinear PDEs are typically difficult to acquire. Therefore, numerical solutions to these kinds of nonlinear PDEs are very important. However, due to their exponential increase in complexity, traditional approaches like finite difference method and finite element method fail in high-dimensional instances. Many fields pay close attention to developments in numerical algorithms for solving high-dimensional PDEs. There are several numerical methods for solving nonlinear high-dimensional partial differential equations here, such as Monte Carlo method\cite{kantas2014sequential,warin2018nesting}, lattice rule\cite{dick2013high} and sparse grid method\cite{shen2010efficient,wang2016sparse}, etc. They exhibit relative adaptability in addressing high-dimensional problems. However, they typically require substantial computational resources, especially in high-dimensional scenarios. Monte Carlo method often demands a large number of sample points, while lattice rule and sparse grid method may require finer grids or adaptive strategies. Moreover, their convergence rates are usually relatively slow, particularly in high-dimensional situations. Achieving the desired accuracy may entail a significant amount of computation.

Recently, deep neural networks (DNNs) have been used to create numerical algorithms which work well at overcoming the curse of dimensionality and successfully solving high-dimensional PDEs\cite{raissi2024forward,beck2021deep,yu2018deep,han2018solving,han2020solving,khoo2021solving,khoo2019switchnet,lu2021deepxde,lye2021iterative,lyu2022mim,raissi2019physics,samaniego2020energy, sirignano2018dgm,wang2020mesh,zang2020weak,zeng2022adaptive,zhu2019physics}. Inspired by Ritz method, deep Ritz method  (DRM) \cite{yu2018deep} is proposed to solve variational problem arising from PDEs. The deep Galerkin method (DGM) is proposed in \cite{sirignano2018dgm} to solve high-dimensional PDEs by approximating the solution with a deep neural network which is trained to satisfy the differential operator, initial condition, and boundary conditions.
 The physics-informed neural networks (PINN) is presented in \cite{raissi2019physics}, where the PDE is embedded into the neural network by utilizing automatic differentiation (AD). 
 Three adaptive techniques to improve the computational performance of DNNs methods for high-dimensional PDEs are presented in \cite{zeng2022adaptive}.
The authors in \cite{khoo2019switchnet} proposed an approach for scattering problems connected with linear PDEs of the Helmholtz type that relies on DNNs to describe the forward and inverse map.
 In \cite{han2018solving}, the deep backward stochastic differential equation (BSDE) method based on the nonlinear Feynman-Kac formula (see e.g.\cite{pardoux1992backward}) is proposed, which is used in \cite{han2020solving} to estimate the solution of eigenvalue problems for semilinear second order differential operators. 

 The PINN and deep BSDE method are two different kinds of numerical frameworks for solving general PDEs. The AD is used to avoid truncation error and the numerical quadrature errors of variational form. Some gradient optimization methods are used to update the neural network so that the loss of the differential equation and boundary condition is reduced. 
 The deep BSDE method treats the BSDE as a stochastic control problem with the gradient of the solution being the policy function and parameterizes the control process of the solution by DNNs. Then it trains the resulting sequence of networks in a global optimization given by the prescribed terminal condition. These methods do not rely on the training data provided by some external algorithms, which can be considered as unsupervised learning methods. One drawback of PINN is the high computational complexity of its loss function, which includes the differential operator in the PDE to be solved. On the other hand, the deep BSDE method does not require the computation of high order derivatives. Moreover, the loss function used by deep BSDE method involves only simple additive calculations, thereby deep BSDE method iterates faster. 
 The deep BSDE method has made high-dimensional problems solvable, which allows us to solve high-dimensional semilinear PDEs in a reasonable amount of time. Recently, there are some works related to deep BSDE method, see \cite{raissi2024forward,beck2021deep,han2020solving,beck2018solving,boussange2023deep,chan2019machine,feng2023deep,fujii2019asymptotic,han2020derivative,nguwi2022deep,nusken2021interpolating,nusken2021solving,pham2021neural,takahashi2022new,vidales2018unbiased,wang2022deep}. Based on the deep BSDE method,  an improved method called FBSNNs is proposed in \cite{raissi2024forward}. The method proceeds by approximating the unknown solution using a fully connected feedforward neural network (FNN) and obtains the required gradient vector by applying AD.
 However, because the nonlinear Feynman-Kac formula is involved in the reformulation procedure, FBSNNs can only handle some specific Cauchy problems without boundary conditions. Then, it is desirable to extend the FBSNNs to other kinds of PDEs and deal with the problems with boundary conditions.
 
 The Navier-Stokes equation is an important equation in fluid dynamics  and the Cahn-Hilliard equation is widely used in multi-phase problems. These equations are difficult to solve due to their complexity. There are many deep learning methods that have been applied to solve these equations in one or two dimensions (see e.g.\cite{mattey2022novel,miyanawala2017efficient,mohan2018deep,raissi2018hidden,wight2020solving,zhu2021local}). However, these methods fail due to excessive complexity when the dimension is more than three. We choose to introduce the FBSNNs presented in \cite{raissi2024forward} to solve these equations. 
 We convert the incompressible Navier-Stokes equations into FBSDEs and then employ FBSNNs to solve them in two or three dimension.
 We develop a suitable numerical method based on the reflection principle to deal with the Neumann boundary condition and handle the Dirichlet boundary condition using the method mentioned in \cite{pardoux1998backward}. We rewrite the Cahn-Hilliard equation into a system of parabolic equations by adding stable terms reasonably, then the numerical solution of the new system is obtained  using the FBSNNs. However, when dealing with mixed boundary condition, the above method should be improved. We utilize an approach which is similar to the method for Dirichlet boundary case, meanwhile we add an extra error item to the final loss function for the Neumann boundary condition.
 The equation can also be solved for periodic boundary condition with techniques involved the periodicity.
 Therefore, we can naturally solve the CHNS equation which is a coupled system of Navier-Stokes and Cahn-Hilliard equations.

 The rest of this article is organized as follows. In Section \ref{sec02}, we introduce FBSDEs, deep BSDE method and FBSNNs method briefly.
 In Section \ref{sec03}, we present the approach to solve incompressible Navier-Stokes equations with different boundary conditions. A methodology is proposed in Section \ref{sec04} to solve Cahn-Hilliard equation with different boundary conditions. In Section \ref{sec05}, the method to solve CHNS system is developed. Numerical experiments are given in Section \ref{sec06} to verify the effectiveness of our methods. Finally,
 conclusions and remarks are given in Section \ref{sec07}.

\section{FBSDEs, deep BSDE method and FBSNNs}
\label{sec02}
\subsection{A brief introduction of FBSDEs}
\label{sec02:sec01}
 The FBSDEs where the randomness in the BSDE driven by a forward stochastic differential equation (SDE), is written in the general form
\begin{equation}\label{eq:BSDE1}
	\left\{
	\begin{aligned}
	   & d X_{s}  = b(s,X_{s}) ds+\sigma(s,X_{s})d W_{s}, \quad s\in[0,T],\\
		 & X_0  = x_0,\\
		 & -d Y_{s} = f(s, X_{s}, Y_{s}, Z_{s}) ds - Z_{s}^T\sigma(s,X_{s})d W_{s}, \quad s\in[0,T],\\
	   & Y_T  = g(X_T),
	\end{aligned} 
	\right.
\end{equation}
where $\{W_s\}_{0\le s\le T}$ is a d-dimensional Brownian motion, $b: [0,T]\times\mathbb{R}^{d}\rightarrow \mathbb{R}^{d}$, $\sigma:[0,T]\times\mathbb{R}^{d}\rightarrow \mathbb{R}^{d \times d} $, $f : [0,T]\times\mathbb{R}^{d}\times \mathbb{R}^{m}\times \mathbb{R}^{d\times m} \rightarrow \mathbb{R}^{m}$ and $g: [0,T]\times\mathbb{R}^{d} \rightarrow \mathbb{R}^{m}$ are all deterministic mappings of time and space, with the fixed $T>0$. We refer to $Z$ as the $control\ process$ according to the stochastic control terminology. In order to guarantee the existence of a unique solution pair $\{(Y_s, Z_s)\}_{0\le s\le T}$ adapted to the augmented natural filtration, the standard well-posedness assumptions of \cite{pardoux1992backward} are required.
Indeed, considering the quasi-linear, parabolic terminal problem
\begin{equation}
	\begin{aligned}
		u_t(t,x)  + \mathcal{L}u(t,x) +f(t,x,u(t,x),\nabla u(t,x)) = 0,\  (t,x)\in[0,T]\times\mathbb{R}^{d},
	\end{aligned}\label{eq:PDE}
\end{equation}
with $u(T,x) = g(x)$ and  $\mathcal{L}$ is the second-order differential operator
\begin{equation}\label{eq:operater}
	\begin{aligned}
		& \mathcal{L}u(t,x) = \frac{1}{2}\sum_{i,j = 1}^{d}a_{i,j}(t,x)\frac{\partial^{2} u(t,x)}{\partial x_{i}\partial x_{j}} + \sum_{i=1}^{d}b_{i}(t,x)\frac{\partial u(t,x)}{\partial x_{i}}, \quad   a_{i,j} = [\sigma\sigma^{T}]_{ij},
	\end{aligned}
\end{equation}
the nonlinear Feynman-Kac formula  indicates that the  solution of \eqref{eq:BSDE1} coincides almost exactly with the solution of \eqref{eq:PDE} (cf., e.g.,
\cite{pardoux1992backward})
\begin{equation}\label{eq:BSDE2}
	\begin{aligned}
		& Y_{s} = u(s,X_{s}), \  Z_{s} =  \nabla u \left(s,X_{s}\right), \  s\in [0,T].
	\end{aligned} 
\end{equation}
As a result, the BSDE formulation offers a stochastic representation to the synchronous solution of a parabolic problem and its gradient, which is a distinct advantage for numerous applications in stochastic control.

\subsection{Deep BSDE method}
\label{sec02:sec02}
Inorder to  review the deep BSDE method proposed in \cite{han2018solving}, we consider the following FBSDEs
\begin{equation}\label{eq:FBSDE}
\left\{
	\begin{aligned}
		 & X_{t}  = x_0+\int_{0}^{t}b(s,X_{s}) ds+\int_{0}^{t}\sigma(s,X_{s})d W_{s}, \\
	   & Y_{t}  = g(X_T)+ \int_{t}^{T}f(s, X_{s}, Y_{s}, Z_{s}) ds - \int_{t}^{T}Z_{s}^T\sigma(s,X_{s})d W_{s},
	\end{aligned} 
 \right.
\end{equation}
which is the integration form  of  \eqref{eq:BSDE1}. Given a partition of the time interval $[0,T]:0 = t_{0}< t_{1}<...<t_{N}=T$, the  Euler-Maruyama scheme is used to discretize for $X_{t}$ and $Y_{t}$
and we have
\begin{equation}\label{eq:FBSDE1}
\left\{
	\begin{aligned}
		 &  X_{t_{n+1}}  =  X_{t_{n}}+b(t_n,X_{t_{n}})\Delta t_{n}+\sigma(t_n,X_{t_{n}})\Delta W_{t_n}, \\
		  & Y_{t_{n+1}}(X_{t_{n+1}})  =  Y_{t_{n}}(X_{t_{n}})- f_{t_n}(X_{t_{n}}, Y_{t_{n}}(X_{t_{n}}), Z_{t_{n}}(X_{t_{n}}))\Delta t_{n}\\ 
            & + Z_{t_{n}}^T(X_{t_{n}})\sigma(t_n,X_{t_{n}})\Delta W_{t_n},
	\end{aligned} 
 \right.
\end{equation}
where $\Delta t_{n} = t_{n+1}-t_{n}= \frac{T}{N}$, $\Delta W_{t_n} = W_{t_{n+1}}-W_{t_{n}}$.
The $Z_{t_{n}}(X_{t_{n}})$ is approximated  by  a FNN  with parameter $\theta_n$ for $n=1, \cdots , N-1$. The initial values $Y_{t_{0}}(X_{t_{0}})$ and $Z_{t_{0}}(X_{t_{0}})$ are treated as trainable parameters in the model. To make $Y_{t_{0}}(X_{t_{0}})$ to approximate $u(t_{0},X_{t_{0}})$, the difference in the matching with a given terminal condition is used to define the expected loss function
\begin{equation}\label{eq:error}
	\begin{aligned}
		& l\left(Y_{t_{0}}(X_{t_{0}}), Z_{t_{0}}(X_{t_{0}}),\theta_{1},..., \theta_{N-1}\right) = \frac{1}{M}\sum_{m=1}^{M}\left\vert (g-Y_{t_{N}})(X_{t_{N},m})\right\vert^2,
	\end{aligned} 
\end{equation}
which represents $M$ different realizations of the underlying Brownian
motion, where the subscript $m$ corresponds to the $m$-th realization of the underlying Brownian
motion. The process is called the deep BSDE method.

\subsection{FBSNNs}
Raissi \cite{raissi2024forward} introduced neural networks called FBSNNs to solve FBSDEs. The unknown solution $u(t, x)$ is approximated by the FNN with the input $(t,x)$ and  the required gradient vector $\nabla u (t, x)$ is attained by applying AD. The parameter $\theta$ of FNN  can be learned by minimizing the loss function given explicitly in equation \eqref{eq:error1}  obtained from discretizing the FBSDEs \eqref{eq:FBSDE} using the Euler-Maruyama scheme 
\begin{equation}\label{eq:FBSDE2}
\left\{
	\begin{aligned}
		& X_{t_{n+1}} = X_{t_{n}}+b(t_n,X_{t_{n}})\Delta t_{n}+\sigma(t_n,X_{t_{n}})\Delta W_{t_n}, \\
		  & \widetilde{Y}_{t_{n+1}}(X_{t_{n+1}})  = Y_{t_{n}}(X_{t_{n}})- f_{t_n}(X_{t_{n}}, Y_{t_{n}}(X_{t_{n}}), Z_{t_{n}}(X_{t_{n}}))\Delta t_{n}\\ 
            & + Z_{t_{n}}^T(X_{t_{n}})\sigma(t_n,X_{t_{n}})\Delta W_{t_n},
	\end{aligned} 
 \right.
\end{equation}
where $Y_{t_{n}}(X_{t_{n}})$ represents the estimated value of $u(t_n,X_{t_{n}})$ given by the FNN  and $\widetilde{Y}_{t_{n+1}}(X_{t_{n+1}})$ is the reference value  corresponding to $Y_{t_{n+1}}(X_{t_{n+1}})$, which is obtained from the calculation in \eqref{eq:FBSDE2}. The loss function is then given by 
\begin{equation}\label{eq:error1}
	\begin{aligned}
		l(\theta) & = \sum_{m=1}^{M} \left[\sum_{n=0}^{N-1}\vert(\widetilde{Y}_{t_{n+1}}- Y_{t_{n+1}})(X_{t_{n+1},m})\vert^2  + \vert (g-Y_{t_{N}})(X_{t_{N},m})\vert^2\right],
	\end{aligned} 
\end{equation}
where the subscripts $m$ is the same meaning as it in \eqref{eq:error}. 

  The deep BSDE method only calculates the value of $u(t_{0},X_{t_{0}})$. This
means that in order to obtain an approximation to $Y_{t_{n}}(X_{t_{n}}) = u(t_{n}, X_{t_{n}})$ at a later time $t_{n} > t_{0}$, we have to retrain the algorithm. Furthermore, the number of the FNNs grows with the number of time steps $N$, which makes training difficult. In this article, we use the FBSNNs. The FNN is expected to be able to  approximate $u(t, x)$ over the entire computational area instead of only one point. That is, we will use multiple initial points to train the FNN. In order to improve the efficiency, the number of Brownian motion trajectories for each initial point is set as $M = 1$.

\section{Deep neural network for solving the incompressible Navier-Stokes equation}
\label{sec03}
\subsection{A class of FBSDEs associated to the incompressible Navier-Stokes equation}
\label{sec03:sec01}
The  Cauchy problem for deterministic backward Navier–Stokes equation for the velocity field of the incompressible and viscous fluid is 
\begin{equation}\label{eq:N-S1}
	\left\{ 
	\begin{aligned}
		 & u_t+\nu\Delta u+(u\cdot\nabla)u+\nabla p+f  = 0,   \  0\le t\le T, \\
		& \nabla\cdot u=0,\ u(T) = g , \\
	\end{aligned}
	\right.
\end{equation}
which is obtained from the classical Navier–Stokes equation via the time-reversing transformation
\begin{equation}\label{eq:timerevers}
	\begin{aligned}
    &(u, p, f)(t,x)   \rightarrow (-u, p, f)(T-t,x), \quad   0\le t\le T.
    \end{aligned}
\end{equation}
Here $u = u(t, x)$ represents the $d$-dimensional velocity field of a fluid, $p = p(t, x)$ is the pressure, $\nu>0$ is the viscosity coefficient, and $f = f (t, x)$ stands for the external force. We now study the  backward Navier–Stokes equation \eqref{eq:N-S1} in $\mathbb{R}^{d}$ with different boundary conditions. 


Then, the PDE \eqref{eq:N-S1} is associated through the nonlinear Feynman-Kac formula to the following FBSDEs
\begin{equation}\label{eq:NSBSDE1}
	\left\{ 
	\begin{aligned}
		 & d X_{s}  =  \sqrt{2\nu} d W_{s}, \quad  s \in [0,T],\\
		 & X_{0}  = x_0,\\
		 & -d Y_{s}  = \left(f(s,X_{s})+\nabla p\left(s,X_{s}\right) + (Y_{s}\cdot\nabla) Y_{s}\right) ds- \sqrt{2\nu}Z_{s}^Td W_{s}, \quad  s \in [0,T],\\
		 & Y_{T}  =  g(X_{T}),\\
	\end{aligned}
	\right.
\end{equation}
where 
\begin{equation*}
	Y_{s} = u\left(s,X_{s}\right), \quad Z_{s} = \nabla u\left(s, X_{s}\right). 
\end{equation*}

\subsection{The algorithm for solving the incompressible Navier-Stokes equation}
\label{sec03:sec02}
\begin{figure}[H]
	\centering 
	\includegraphics[scale=0.5]{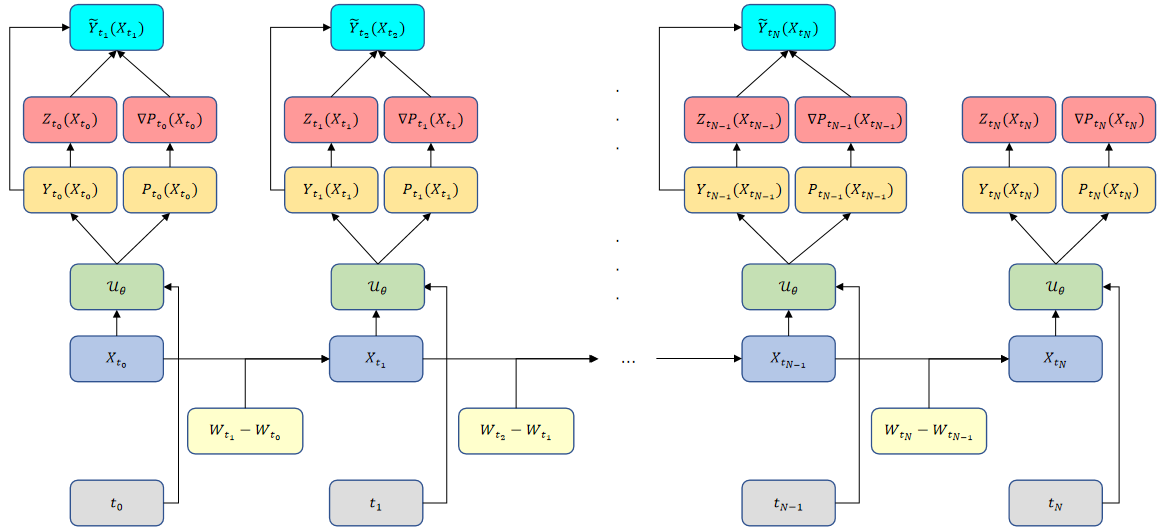}	
	\caption{Illustration of  Algorithm \ref{alg:Algorithm 1} for solving the incompressible Navier-Stokes equation}
	\label{fig:fig1} 
\end{figure}

Given a partition of $[0,T]: 0 = t_{0}< t_{1}<...<t_{N}=T$,  we consider the Euler-Maruyama scheme with $n = 0,..., N-1$ for FBSDEs \eqref{eq:NSBSDE1}
\begin{equation}\label{eq:NSEuler1}
	\left\{ 
	\begin{aligned}
		 & X_{t_{n+1}}  = X_{t_{n}} + \sqrt{2\nu}\Delta W_{t_n},\\
		 & \widetilde{Y}_{t_{n+1}}(X_{t_{n+1}})  = Y_{t_{n}}(X_{t_{n}})-(f_{t_{n}}+ \nabla P_{t_{n}}+ (Y_{t_{n}}\cdot \nabla)Y_{t_{n}})(X_{t_{n}})\Delta t_{n}\\  
           & + \sqrt{2\nu} Z_{t_{n}}^T(X_{t_{n}}) \Delta W_{t_n},
	\end{aligned}
	\right.
\end{equation}
where $\Delta t_{n} = t_{n+1}-t_{n}=\frac{T}{N}$ and $\Delta W_{t_n} = W_{t_{n+1}}-W_{t_{n}}$. The 
$(Y_{t_{n}},P_{t_{n}})^T$ represents the estimated value of 
$(u,p)^T$ at time $t_n$ given by the FNN, respectively. The $\widetilde{Y}_{t_{n+1}}(X_{t_{n+1}})$ is the reference value of $Y_{t_{n+1}}(X_{t_{n+1}})$, which is obtained from the calculation in the second equation in \eqref{eq:NSEuler1}. We utilize $K$ different initial sampling points for training the FNN. The algorithm  of the proposed scheme is summarized in Algorithm \ref{alg:Algorithm 1}. Illustration of the Algorithm \ref{alg:Algorithm 1} for solving the incompressible Navier-Stokes equation is shown in Figure \ref{fig:fig1}.
\begin{algorithm}
	\caption{Algorithm  for  the
		incompressible Navier-Stokes equation}
	\label{alg:Algorithm 1}
	\begin{algorithmic}[1]
		\Require
		 Number of initial sampling points $K$,   terminal time $T$, number of time intervals $N$, viscosity coefficient $\nu$, maximum number of training iterations $M_{iter}$.
		\Ensure
		The optimal FNN $\mathcal{U}_\theta$.
		\State  Initialize the FNN $\mathcal{U}_\theta$;
            \State  Select initial sampling points $x_0$ by uniform distribution;
		\State Generate independent $d$-dimensional standard Brownian motions $W_{t_n}(n = 0,...,N)$;
		\State Compute $X_{t_{n+1}}$ according to \eqref{eq:NSEuler1} for $n = 0,..., N-1$;
            \State According to \eqref{eq:NSEuler1}, use the FNN $\mathcal{U}_\theta$ with AD to calculate $\widetilde{Y}_{t_{n+1}}(X_{t_{n+1}})$ for $n = 0, \ldots, N-1$;
		\State Minimize the loss function by the Adam algorithm
    \begin{equation}\label{eq:lossNSBSDE1}
	\begin{aligned}
        l(\theta) & =  \frac{1}{K}\sum_{k=1}^{K}\left[\frac{1}{N}\sum_{n=0}^{N-1}\vert(\widetilde{Y}_{t_{n+1}}-Y_{t_{n+1}})(X_{t_{n+1},k})\vert^2 +\alpha_1\vert (g-Y_{t_{N}})(X_{t_{N},k})\vert^{2} \right.\\ 
        & \left.+\frac{\alpha_2}{N+1}\sum_{n=0}^{N}\vert\nabla\cdot Y_{t_{n}}(X_{t_{n},k})\vert^2\right],\\
            \end{aligned}
    \end{equation}
    where $\alpha_i,i=1,2$ are the weights of the components of the loss function. The subscript $k$ corresponds to the $k$-th initial sampling point;
		\State Repeat procedures 2 to 6 until the maximum number  of training iterations $M_{iter}$ is reached.
	\end{algorithmic}
\end{algorithm}


\subsection{The algorithm for solving the incompressible Navier-Stokes equation with Dirichlet boundary condition}
\label{sec03:sec03}
For the backward Navier–Stokes equation \eqref{eq:N-S1} in $\Omega \subset \mathbb{R}^d$ with the Dirichlet boundary condition 
\begin{equation}
	\begin{aligned}
		& u(t,x) = h(t,x),\quad (t,x)\in[0,T]\times\partial\Omega,
	\end{aligned}\label{eq:dboundary}
\end{equation} 
the corresponding FBSDEs can be rewritten as the following form according to \cite{pardoux1998backward} by the nonlinear Feynman-Kac formula
\begin{equation}
	\left\{ 
	\begin{aligned}
		 & X_t  =x_0+\int_{0}^{t}\sqrt{2\nu}dW_s,\\
          & Y_t  =\Phi(T\Lambda\tau,X_{T\Lambda\tau})+\int_{t\Lambda\tau}^{T\Lambda\tau}(f(s,X_s)+\nabla p\left(s,X_s\right)\\
		& +(Y_s\cdot\nabla)Y_s)ds
		-\int_{t\Lambda\tau}^{T\Lambda\tau}\sqrt{2\nu}Z_s^TdW_s,
	\end{aligned}\label{eq:dboundaryFBSDEs}
	\right. 
\end{equation}
where $a\Lambda b = \min\{a,b\}$, the stopping time $\tau = \inf\{ t>0:X_{t}\notin \Omega\}$ be the first time that the process $X_{t}$ exits $\Omega$ and 
\begin{equation}\label{eq:BCdef}
	\begin{aligned}
		\Phi(T\Lambda\tau,X_{T\Lambda\tau})=\left\{
		\begin{array}{ll}
			g(X_T),\quad\tau> T,\ X_T\in \Omega, \\
			h(\tau,X_\tau),\quad\tau\le T,\ X_\tau\in\partial\Omega.\\
		\end{array}
		\right.
	\end{aligned}
\end{equation}
Through the Euler scheme, the discrete formulation of the FBSDEs \eqref{eq:dboundaryFBSDEs} can be obtained accordingly
\begin{equation}\label{eq:NSEuler2}
	\left\{ 
	\begin{aligned}
		& X_{t_{n+1}}  = X_{t_{n}} + \sqrt{2\nu}\Delta W_{t_n},\\
		& \widetilde{Y}_{t_{n+1}\Lambda\tau}(X_{t_{n+1}\Lambda\tau})  =Y_{t_{n}\Lambda\tau}(X_{t_{n}\Lambda\tau})-\mathbbm{1}_{(0,\tau)}(t_{n})[(f_{t_{n}}+\nabla P_{t_{n}}\\
            & +(Y_{t_{n}}\cdot\nabla)Y_{t_{n}})(X_{t_{n}})\Delta t_{n}-\sqrt{2\nu}Z_{t_{n}}^T(X_{t_{n}})\Delta W_{t_n}],
	\end{aligned}
	\right.
\end{equation}
where $\mathbbm{1}_{(0,\tau)}(t_{n})=1,t_{n}\in[0,\tau)$. It should be noted that we will calculate the stop time $\tau$ after the iteration of $X_{t_{n+1}}$ is completed. Supposing $\tau=t_{n+1}$ when $\tau\le T$, we let $X_{\tau} = X_{t_{n+1}}$ on $\partial\Omega$ and update $\Delta W_{t_{n}}=(X_{t_{n+1}}-X_{t_{n}})/\sqrt{2\nu}$ to satisfy \eqref{eq:NSEuler2}.  
The algorithm of the proposed scheme is similar as Algorithm \ref{alg:Algorithm 1}.

\subsection{The algorithm for solving the incompressible Navier-Stokes equation with Neumann boundary condition}
\label{sec03:sec04}
We consider the backward Navier–Stokes equation \eqref{eq:N-S1} with the Neumann boundary condition 
\begin{equation}
	\begin{aligned}
		& \frac{\partial u(t,x)}{\partial \mathbf{n}}=q(t,x),\quad (t,x)\in[0,T]\times\partial\Omega,
	\end{aligned}\label{eq:nboundary}
\end{equation} 
where $\mathbf{n}$ is the unit normal vector at $\partial \Omega$ pointing outward of $\Omega$. Supposing $x\in\partial\Omega, x+\Delta x\notin\Omega$ and $x-\Delta x\in\Omega$, where $x+\Delta x$ and $x-\Delta x$ are symmetric to the boundary $\partial \Omega$. Then we have
\begin{equation}
	\begin{aligned}
		& u(t,x-\Delta x)-u(t,x+\Delta x) \approx -2q(t,x)\vert\Delta x\vert,\quad (t,x)\in[0,T]\times\partial\Omega.
	\end{aligned}\label{eq:reflect}
\end{equation} 
If $X_{t_{n+1}}\in\Omega$, let $X'_{t_{n+1}}=X_{t_{n+1}}$, and if $X_{t_{n+1}}\notin\Omega$, let $X'_{t_{n+1}}\in\Omega$ is the symmetric point of $X_{t_{n+1}}$ to the boundary $\partial \Omega$. The $X''_{t_{n+1}}$ is used to denote the intersection of the line segment $\overline{X_{t_{n+1}}X'_{t_{n+1}}}$ and $\partial \Omega$. Therefore, the discretization can be rewritten similarly as
\begin{equation}
	\left\{ 
	\begin{aligned}
	   & X_{t_{n+1}}  = X_{t_{n}} + \sqrt{2\nu}\Delta W_{t_n},\\
		 & \Delta Y_{t_{n+1}}  = q(t_{n+1},X''_{t_{n+1}})\vert\Delta X_{t_{n+1}}\vert,\\
	   & X_{t_{n+1}}  = X'_{t_{n+1}},\\
	   &	\widetilde{Y}_{t_{n+1}}(X_{t_{n+1}})  = Y_{t_{n}}(X_{t_{n}})-(f_{t_{n}}+\nabla P_{t_{n}}+(Y_{t_{n}}\cdot\nabla)Y_{t_{n}})(X_{t_{n}})\Delta t_{n} \\
         & +\sqrt{2\nu}Z_{t_{n}}^T(X_{t_{n}})\Delta W_{t_n}-\Delta Y_{t_{n+1}},
	\end{aligned}\label{eq:NSEuler3}
	\right.
\end{equation}
where $\Delta X_{t_{n+1}}= X_{t_{n+1}}-X'_{t_{n+1}}$. 
The algorithm of the proposed scheme is similar as Algorithm \ref{alg:Algorithm 1}.
    \begin{remark}
There are some similar works in \cite{boussange2023deep,han2020derivative} for dealing with the Neumann boundary conditions.
If $X_{t_{n+1}}\notin\Omega$ during the iterative process, the authors \cite{boussange2023deep, han2020derivative} choose to reflect $X_{t_{n+1}}$  on the boundary $\partial\Omega$, which allows them to deal with the homogeneous Neumann conditions. In contrast, our method can deal with non-homogeneous Neumann boundary conditions.
    \end{remark}

\section{Deep neural network for solving the Cahn-Hilliard equation}
\label{sec04}
\subsection{Rewrite Cahn-Hilliard equation into  a parabolic PDE system}
\label{sec04:sec01}
We consider the following Cahn-Hilliard equation, which has fourth order derivatives,  
\begin{equation}\label{eq:CH}
	\left\{
	\begin{aligned}
		&  \phi_t -L_d\Delta \mu + f= 0, \quad t\ge 0, \\
		& \mu + \gamma^2\Delta\phi+\phi-\phi^3 = 0, \quad t\ge 0,\\
		& \phi(0) = -g, 
	\end{aligned}
	\right.
\end{equation}
 where $\phi= \phi(t, x)$ is the unknown, e.g., the concentration of the fluid, $\mu = \mu(t, x)$ is a function of $\phi$, e.g.,the chemical potential, $L_{d} > 0$ is the diffusion coefficient and $\gamma> 0$ is the model parameter. 
A first order stabilized scheme \cite{Shen2010stab} for the Cahn-Hilliard  equation \eqref{eq:CH} reads as 
\begin{equation}\label{eq:CH2}
	\left\{
	\begin{aligned}
		 & \phi^{n+1} -\phi^{n} -\Delta t L_d\Delta\mu^{n} + \Delta t f^{n}  = 0, \\
		 & \mu^{n+1} + \gamma^2\Delta \phi^{n}+\phi^{n}-(\phi^{n})^3-\frac{S}{L_d}(\phi^{n+1} -\phi^{n})  = 0,\\ 
	\end{aligned}
	\right.
\end{equation}
where $S$ is a suitable stabilized parameter. 
It is easy to derive the following equation
\begin{equation}\label{eq:CH3}
	\begin{aligned}
		& \frac{ \mu^{n+1}-\mu^{n}}{\Delta t} + \frac{\gamma^2\Delta \phi^{n}}{\Delta t} +\frac{\phi^{n}-(\phi^{n})^3-\frac{S}{L_d}(\phi^{n+1} -\phi^{n})}{\Delta t} =\frac{-\mu^{n}}{\Delta t}. 
	\end{aligned}
\end{equation}
The first equation of \eqref{eq:CH2}  and \eqref{eq:CH3} can be regarded as the discretization of the following modified Cahn-Hilliard equation in $[0,T]$
\begin{equation}\label{eq:CH4}
	\left\{
	\begin{aligned}
		&  \phi_t - L_d\Delta\mu + f= 0, \\
		& \mu_t + \frac{\gamma^2}{\delta}\Delta \phi-S\Delta \mu + \frac{S}{L_d}f+\frac{1}{\delta}\left(\mu + \phi-\phi^3\right) = 0, \\
		& \phi(0) = -g, 
	\end{aligned}
	\right.
\end{equation}
where $\delta = O (\Delta t)$.
By reversing the time and defining $$(\phi, \mu, f)(t,x) \rightarrow (-\phi, -\mu,   f)(T-t,x),\quad   0\le t\le T,$$
the $(\phi, \mu)$ satisfies the following backward Cahn-Hilliard equation in  $[0,T]$
\begin{equation}\label{eq:CH5}
	\left\{
	\begin{aligned}
		&  \phi_t + L_d\Delta\mu + f= 0, \\
		& \mu_t - \frac{\gamma^2}{\delta}\Delta\phi+S\Delta\mu+\frac{S}{L_d}f-\frac{1}{\delta}(\mu + \phi-\phi^3) = 0,\\
		& \phi(T) = g.  
	\end{aligned}
	\right.
\end{equation}
In order to satisfy the nonlinear Feynman-Kac formula and utilize the FBSNNs, we treat the backward Cahn-Hilliard  equation \eqref{eq:CH5} as a semilinear parabolic differential equation
\begin{equation}\label{eq:CHMatrix}
	\begin{aligned}
		\psi_{t} + A\Delta\psi + F = 0,
	\end{aligned}
\end{equation}
with $\psi = (\phi, \mu)^{T}$, $A =\left(\begin{array}{ll} 0 &  L_{d}  \\ -\frac{\gamma^2}{\delta}& S \end{array} \right)$ and $F=\left(f,\frac{S}{L_d}f-\frac{1}{\delta}(\mu + \phi-\phi^3)\right)^{T}$. Supposing 
$\lambda_1$ and $\lambda_2$ are two different eigenvalues of the coefficient matrix $A$,
 then the coefficient matrix $A$ can be diagonalized by $R$ and $R^{-1}$ so that $D = \mbox{diag}(\lambda_1,\lambda_2) = R^{-1}AR$, where $R$ is a matrix of eigenvectors. The system \eqref{eq:CHMatrix} becomes
\begin{equation}\label{eq:CHMatrixD}
	\begin{aligned}
		\hat{\psi}_{t} + D\Delta\hat{\psi} + \hat{F} = 0,
	\end{aligned}
\end{equation}
where  
$\hat{\psi} = R^{-1}\psi=(\hat{\phi},\hat{\mu})^{T}$, $\hat{F} = R^{-1}F=(\hat{F}^{\hat{\phi}},\hat{F}^{\hat{\mu}})^{T}$ and $S$ is chosen so that $ S> 2\gamma\sqrt{\frac{ L_{d}}{\delta}}$ for $\lambda_{1,2} > 0$, where $\lambda_{1,2} = \frac{S \pm \sqrt{S^{2}-\frac{4\gamma^2 L_d}{\delta}}}{2}$. 

Therefore, the system \eqref{eq:CHMatrixD} is decomposed into two independent PDEs and the corresponding FBSDEs can be obtained as follows
\begin{equation}\label{eq:CHBSDE1}
	\left\{
	\begin{aligned}
		& dX_{s}^{\hat{\phi}}=\sqrt{2\lambda_1}dW_{s}^{\hat{\phi}}, \quad  s \in [0,T],\\
		& dX_{s}^{\hat{\mu}}=\sqrt{2\lambda_2}dW_{s}^{\hat{\mu}}, \quad  s \in [0,T],\\
		& X_{0}^{\hat{\phi}} = x_0,\ X_{0}^{\hat{\mu}} = x_0,\\
		& -dY_{s}^{\hat{\phi}}=\hat{F}_{s}^{\hat{\phi}}ds-\sqrt{2\lambda_1}(Z_{s}^{\hat{\phi}})^TdW_{s}^{\hat{\phi}}, \quad  s \in [0,T],\\
		&  -dY_{s}^{\hat{\mu}}=\hat{F}_{s}^{\hat{\mu}}ds-\sqrt{2\lambda_2}(Z_{s}^{\hat{\mu}})^TdW_{s}^{\hat{\mu}}, \quad  s \in [0,T],\\
		& \phi(T) = g,
	\end{aligned}
	\right.
\end{equation}
where 
\begin{equation*}
	Y_{s}^{\hat{\phi}} = \hat{\phi}\left(s,X_{s}^{\hat{\phi}}\right), Z_{s}^{\hat{\phi}} = \nabla \hat{\phi}\left(s, X_{s}^{\hat{\phi}}\right), \hat{F}_{s}^{\hat{\phi}} = \hat{F}^{\hat{\phi}}\left(s,X_{s}^{\hat{\phi}},\phi(s,X_{s}^{\hat{\phi}}),\mu(s,X_{s}^{\hat{\phi}})\right), 
\end{equation*}
and
\begin{equation*}
	Y_{s}^{\hat{\mu}} = \hat{\mu}\left(s, X_{s}^{\hat{\mu}}\right), Z_{s}^{\hat{\mu}} = \nabla \hat{\mu}\left(s,X_{s}^{\hat{\mu}}\right),
 \hat{F}_{s}^{\hat{\mu}} = \hat{F}^{\hat{\mu}}\left(s,X_{s}^{\hat{\mu}},\phi(s,X_{s}^{\hat{\mu}}),\mu(s,X_{s}^{\hat{\mu}})\right).
\end{equation*}
The $\{X_{s}^{\hat{\phi}}\}_{0\le s\le T}$ and $\{X_{s}^{\hat{\mu}}\}_{0\le s\le T}$ are the forward stochastic processes corresponding to $\hat{\phi}$ and $\hat{\mu}$ respectively, which are constrained by $x_0$ at the initial time.

\subsection{The algorithm for solving the Cahn-Hilliard equation}
\label{sec4:sec02}
Given a partition of $[0,T]: 0 = t_{0}< t_{1}<...<t_{N}=T$, we consider the simple Euler scheme for the FBSDEs \eqref{eq:CHBSDE1} with $n = 0,..., N-1$
\begin{equation}\label{eq:CHEuler0}
      \left\{
	\begin{aligned}
		& X_{t_{n+1}}^{\hat{\phi}}=X_{t_{n}}^{\hat{\phi}}+\sqrt{2\lambda_1}\Delta W_{t_{n}}^{\hat{\phi}},\\
		& X_{t_{n+1}}^{\hat{\mu}}=X_{t_{n}}^{\hat{\mu}}+\sqrt{2\lambda_2}\Delta W_{t_{n}}^{\hat{\mu}}.\\
		&\widetilde{Y}_{t_{n+1}}^{\hat{\phi}}(X_{t_{n+1}}^{\hat{\phi}})=Y_{t_{n}}^{\hat{\phi}}(X_{t_{n}}^{\hat{\phi}})-\hat{F}_{t_{n}}^{\hat{\phi}}(X_{t_{n}}^{\hat{\phi}},Y_{t_{n}}^{\phi}(X_{t_{n}}^{\hat{\phi}}), Y_{t_{n}}^{\mu}(X_{t_{n}}^{\hat{\phi}}))\Delta t_{n}\\
            &+\sqrt{2\lambda_1}(Z_{t_{n}}^{\hat{\phi}}(X_{t_{n}}^{\hat{\phi}}))^T\Delta W_{t_{n}}^{\hat{\phi}},\\
		&\widetilde{Y}_{t_{n+1}}^{\hat{\mu}}(X_{t_{n+1}}^{\hat{\mu}})=Y_{t_{n}}^{\hat{\mu}}(X_{t_{n}}^{\hat{\mu}})-\hat{F}_{t_{n}}^{\hat{\mu}}(X_{t_{n}}^{\hat{\mu}},Y_{t_{n}}^{\phi}(X_{t_{n}}^{\hat{\mu}}), Y_{t_{n}}^{\mu}(X_{t_{n}}^{\hat{\mu}}))\Delta t_{n}\\
            &+\sqrt{2\lambda_2}(Z_{t_{n}}^{\hat{\mu}}(X_{t_{n}}^{\hat{\mu}}))^T\Delta W_{t_{n}}^{\hat{\mu}},
      \end{aligned}
       \right.
\end{equation}		
where  $\Delta t_{n}=t_{n+1}-t_{n}= \frac{T}{N} =\delta ,\  \Delta W_{t_{n}}^{\hat{\phi}} = W_{t_{n+1}}^{\hat{\phi}}-W_{t_{n}}^{\hat{\phi}}$ and $\Delta W_{t_{n}}^{\hat{\mu}} = W_{t_{n+1}}^{\hat{\mu}}-W_{t_{n}}^{\hat{\mu}}$. The $(Y_{t_{n}}^{\hat{\phi}}(X_{t_{n}}^{\hat{\phi}}), Y_{t_{n}}^{\hat{\mu}}(X_{t_{n}}^{\hat{\mu}}))^T$ represents the estimated value of $(\hat{\phi}(t_n,X_{t_{n}}^{\hat{\phi}}),\hat{\mu}(t_n,X_{t_{n}}^{\hat{\mu}}))^T$. The $(\widetilde{Y}_{t_{n+1}}^{\hat{\phi}}(X_{t_{n+1}}^{\hat{\phi}}),\widetilde{Y}_{t_{n+1}}^{\hat{\mu}}(X_{t_{n+1}}^{\hat{\mu}}))^T$ is the reference value of $(Y_{t_{n+1}}^{\hat{\phi}}(X_{t_{n+1}}^{\hat{\phi}}),$ 
$ Y_{t_{n+1}}^{\hat{\mu}}(X_{t_{n+1}}^{\hat{\mu}}))^T$,  which is obtained from the last two equations in \eqref{eq:CHEuler0}.

\begin{figure}[H]
	\centering 
	\includegraphics[scale=0.28]{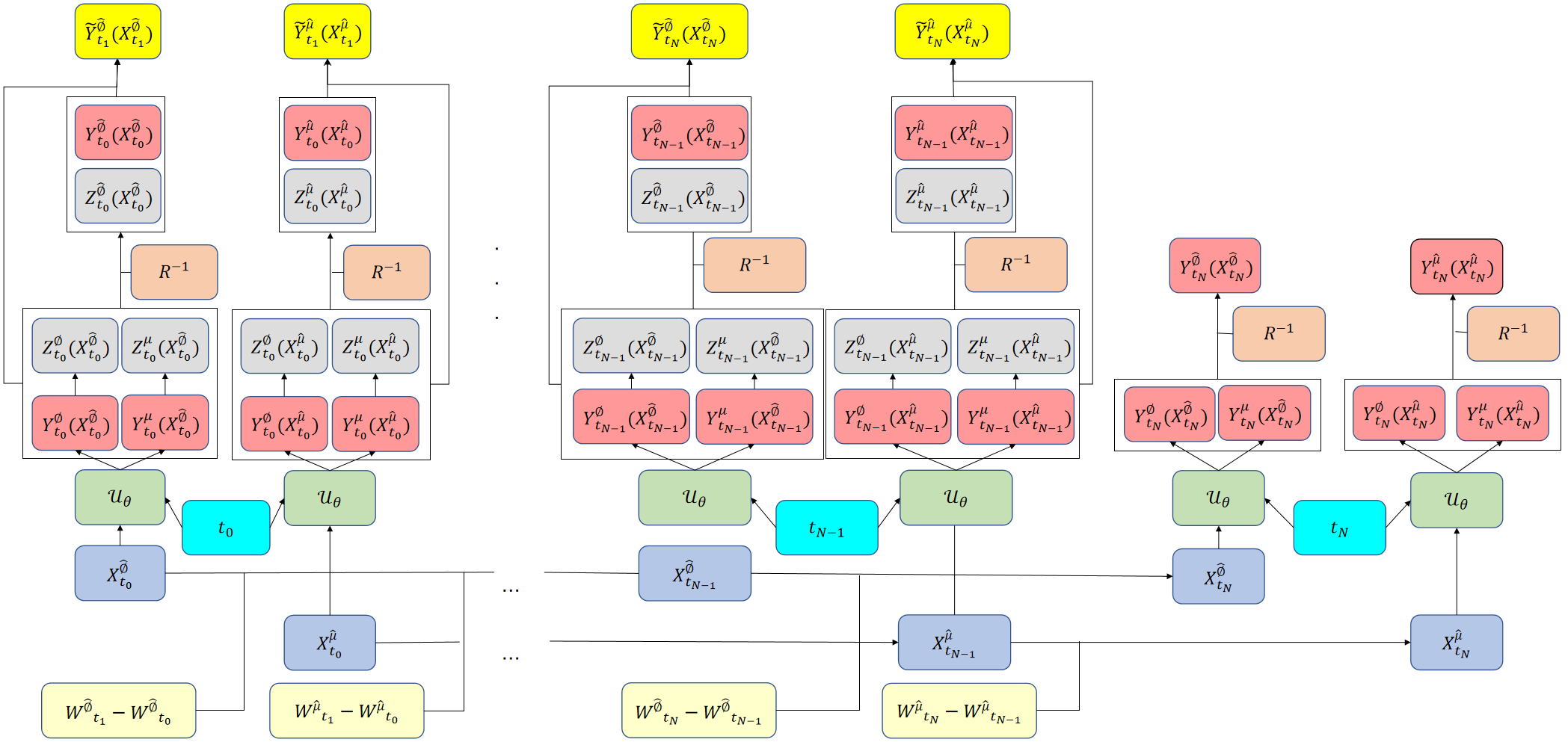}	
	\caption{Illustration of Algorithm \ref{alg:Algorithm 2} for solving the Cahn-Hilliard equation}
	\label{fig:fig2} 
\end{figure}

The $(Y_{t_{n}}^{\phi}, Y_{t_{n}}^{\mu})^T$ represents the estimated value of $(\phi,\mu)^T$ at time $t_n$  given by the FNN. 
Due to the diagonalization, we have 
\begin{equation}\label{eq:diag}
\left\{
	\begin{aligned}
		\begin{aligned}
			&Y_{t_{n}}^{\hat{\phi}}(X_{t_{n}}^{\hat{\phi}}) = \left[R^{-1}(Y_{t_{n}}^{\phi},Y_{t_{n}}^{\mu})^{T}(X_{t_{n}}^{\hat{\phi}})\right]_1,  \\
                &Y_{t_{n}}^{\hat{\mu}}(X_{t_{n}}^{\hat{\mu}}) = \left[R^{-1}(Y_{t_{n}}^{\phi},Y_{t_{n}}^{\mu})^{T}(X_{t_{n}}^{\hat{\mu}})\right]_2,
		\end{aligned}
	\end{aligned}
 \right.
\end{equation}
 where subscript $1$ or $2$  is used to represent the $1$-th  or $2$-th component of the vector. We utilize $K$ different initial sampling points for training the FNN.
The algorithm  of the proposed scheme is summarized as Algorithm \ref{alg:Algorithm 2}. Illustration of the Algorithm \ref{alg:Algorithm 2} for solving the Cahn-Hilliard equation is shown in Figure \ref{fig:fig2}. 
\begin{algorithm}
	\caption{ Algorithm  for  the
		Cahn-Hilliard equation}
	\label{alg:Algorithm 2}
	\begin{algorithmic}[1]
		\Require
		 Number of initial sampling points $K$,  terminal time $T$, number of time intervals $N$,  diffusion coefficient $L_d$, parameters $\gamma$, $\delta$ and $S$, matrix $R^{-1}$,
		eigenvalues $\lambda_{1,2}$, maximum number of training iterations $M_{iter}$.
		\Ensure
		The optimal FNN $\mathcal{U}_\theta$.
		\State Initialize the FNN $\mathcal{U}_\theta$;
            \State  Select initial sampling points $x_0$ by uniform distribution;
		\State Generate independent $d$-dimensional standard Brownian motions $W_{t_n}^{\hat{\phi}}$ and $W_{t_n}^{\hat{\mu}}(n = 0,...,N)$;
		\State  Compute $X_{t_{n+1}}^{\hat{\phi}}$ and $X_{t_{n+1}}^{\hat{\mu}}$ according to \eqref{eq:CHEuler0} for $n = 0,..., N-1$;
            \State According to \eqref{eq:CHEuler0} and \eqref{eq:diag}, use the FNN $\mathcal{U}_\theta$ with AD to calculate  $\widetilde{Y}_{t_{n+1}}^{\hat{\phi}}(X_{t_{n+1}}^{\hat{\phi}})$ and $\widetilde{Y}_{t_{n+1}}^{\hat{\mu}}(X_{t_{n+1}}^{\hat{\mu}})$ for $n = 0, \ldots, N-1$;
		\State Minimize the loss function by the Adam algorithm
		\begin{equation}\label{eq:lossCHBSDE1}
    	\begin{aligned}
    		l(\theta) = & \frac{1}{K}\sum_{k=1}^{K}\left[\frac{1}{N}\sum_{n=0}^{N-1}\left(\vert (\widetilde{Y}_{t_{n+1}}^{\hat{\phi}}-Y_{t_{n+1}}^{\hat{\phi}})(X_{t_{n+1},k}^{\hat{\phi}})\vert^{2}\right.\right.\\
    		& \left.\left.+\vert (\widetilde{Y}_{t_{n+1}}^{\hat{\mu}}-Y_{t_{n+1}}^{\hat{\mu}})(X_{t_{n+1},k}^{\hat{\mu}})\vert^{2}\right)+\alpha_1\vert (g-Y_{t_N}^{\phi})(X_{t_{N},k}^{\hat{\phi}\cup\hat{\mu}})\vert^{2}\right],
    	\end{aligned}
        \end{equation}
        where $X_{t_{N},k}^{\hat{\phi}\cup\hat{\mu}} = X_{t_{N},k}^{\hat{\phi}}\cup X_{t_{N},k}^{\hat{\mu}}$ and $\alpha_1$ is the weight of the terminal condition. The subscript $k$ corresponds to the $k$-th initial sampling point;
		\State Repeat procedures 2 to 6 until the maximum number of training iterations $M_{iter}$ is reached.
	\end{algorithmic}
\end{algorithm}
The Algorithm \ref{alg:Algorithm 2} shows that we only need to compute first-order derivatives during training, which causes the training time to increase linearly with the dimension $d$. This makes our method capable of efficiently solving high-dimensional Cahn-Hilliard equations, which can be observed in the numerical experiments of solving the high-dimensional Cahn-Hilliard equations in Section \ref{sec06:sec02}.

\subsection{The algorithm for solving the Cahn-Hilliard equation with mixed boundary condition}
\label{sec04:sec03}
We consider the  Cahn-Hilliard equation \eqref{eq:CH} in $\Omega \subset \mathbb{R}^d$ with the mixed condition
\begin{equation}\label{eq:dnboundary1}
	\begin{aligned}
		& \left\{
		\begin{aligned}
			&\phi(t,x) = h(t,x), \\
			&\frac{\partial \mu(t,x)}{\partial \mathbf{n}} = q(t,x),
		\end{aligned}
		\right. \quad(t,x)\in[0,T]\times\partial\Omega,
	\end{aligned}
\end{equation}
where $\mathbf{n}$ is the unit normal vector at $\partial\Omega$ pointing outward of $\Omega$. The method described in Section \ref{sec03:sec03} is used to deal with the Dirichlet boundary condition. For the  Neumann boundary condition, it is noted that
\begin{equation}\label{eq:diag2}
	\begin{aligned}
		\begin{aligned}
			\frac{\partial \hat{\mu}}{\partial \mathbf{n}}(t,x) = \left[R^{-1}\left(\frac{\partial \phi}{\partial \mathbf{n}},\frac{\partial \mu}{\partial \mathbf{n}}\right)^T(t,x)\right]_2
		\end{aligned}
	\end{aligned}
\end{equation}
where $\frac{\partial \phi}{\partial \mathbf{n}}(t,x)$ is given by the FNN with AD and subscript 2 represents the second component. Therefore, the method described in Section \ref{sec03:sec04} can be used to deal with the Neumann boundary condition. It is shown in Section \ref{sec06} that our method performs well numerically. The algorithm of the proposed scheme is similar as Algorithm \ref{alg:Algorithm 2}.

\subsection{The algorithm for solving the Cahn-Hilliard equation with periodic boundary condition}
\label{sec04:sec04}
\begin{figure}[H]
	\centering 
	\includegraphics[scale=0.45]{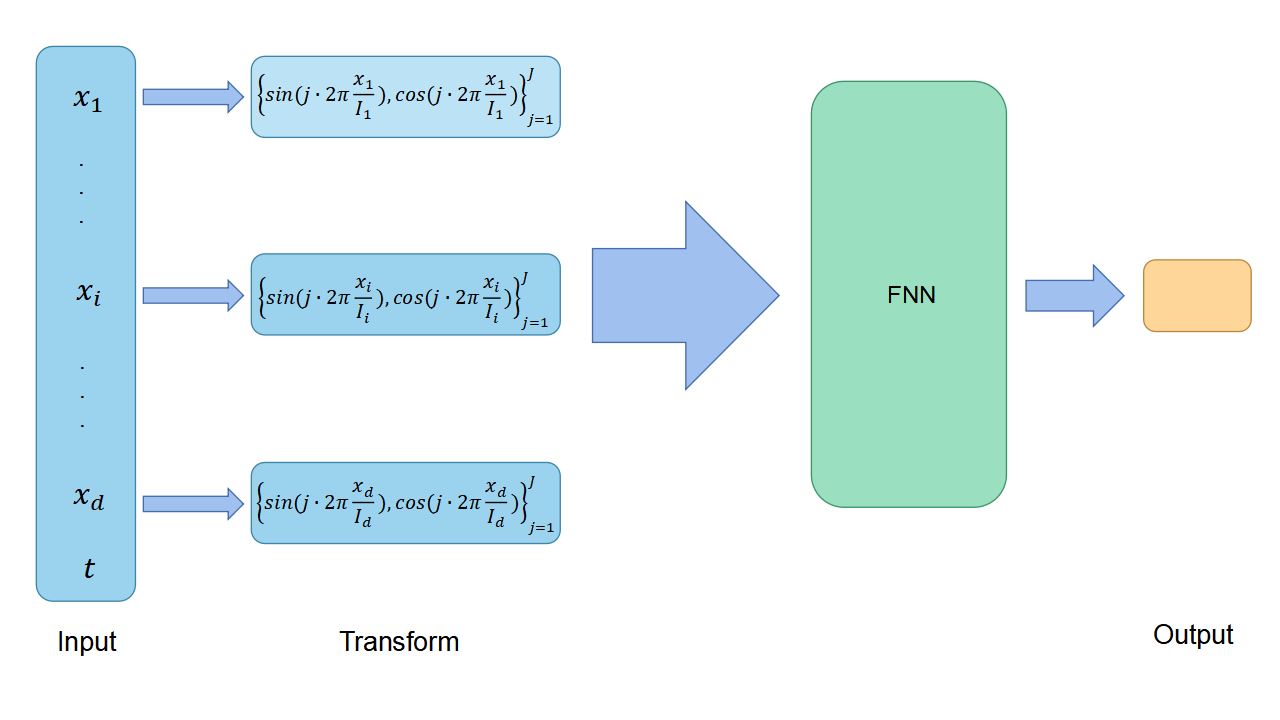}	
	\caption{Network structure for periodic boundary condition}
	\label{fig:fig11} 
\end{figure}

We consider the Cahn-Hilliard equation \eqref{eq:CH} with the periodic boundary condition, 
\begin{equation}\label{eq:pboundary1}
	\begin{aligned}
		& \left\{
		\begin{aligned}
			&\phi(t,x_1,\cdots,x_i+I_i,\cdots,x_d) = \phi(t,x_1,\cdots,x_i,\cdots,x_d),\\
			&\mu(t,x_1,\cdots,x_i+I_i,\cdots,x_d) = 
			\mu(t,x_1,\cdots,x_i,\cdots,x_d), 
		\end{aligned}
		\right. \quad i=1,\cdots,d,
	\end{aligned}
\end{equation}
where $I_i$ is the period along the $i$-th direction. To satisfy the condition \eqref{eq:pboundary1}, as in \cite{han2020solving}, we transform the input vector $x=(x_1,\cdots,x_d)$ into a fixed trigonometric basis before applying
the FNN. The component $x_i$
in $x$ is mapped as follows
\begin{equation}\label{eq:transform}
	\begin{aligned}
		&x_i\rightarrow\left\{\sin\left(j\cdot2\pi \frac{x_i}{I_i}\right),\cos\left(j\cdot2\pi\frac{x_i}{I_i}\right)\right\}_{j=1}^J, \ i=1,\cdots,d,
	\end{aligned}
\end{equation}

where $J$ is the order of the trigonometric basis. The network structure for the periodic boundary condition \eqref{eq:pboundary1} is shown in Figure \ref{fig:fig11}. The algorithm of the proposed scheme is similar as Algorithm \ref{alg:Algorithm 2}.

\section{Deep neural network for solving Cahn-Hilliard-Navier-Stokes system}
\label{sec05}
We now solve the coupled Cahn-Hilliard-Navier-Stokes equation in domain $\mathbb{R}^{d}$. According to Section \ref{sec03:sec01} and Section \ref{sec04:sec01},
after time-reversing, the modified CHNS system is 
\begin{equation}\label{eq:CHNS2}
	\left\{
	\begin{aligned}
		&  u_t+\nu\Delta u+(u\cdot\nabla)u+\nabla p+C\phi\nabla\mu+f_1=0,\\
		&  \phi_t +u\cdot\nabla\phi + L_d\Delta\mu + f_2= 0,  \\
		& \mu_t - \frac{\gamma^2}{\delta}\Delta\phi+S\Delta\mu+\frac{S}{L_d}(u\cdot\nabla\phi + f_2)-\frac{1}{\delta}(\mu + \phi-\phi^3) = 0,\\
            & \nabla \cdot u=0,\\
		& u(T) = g_{u},\quad  \phi(T) = g_{\phi},
	\end{aligned}
	\right.
\end{equation}
 where $C$ denotes a parameter, e.g., the strength of the capillary force comparing with the Newtonian fluid stress. 
 
 Similarly, we have the corresponding FBSDEs of \eqref{eq:CHNS2} by diagonalizing and  using the nonlinear Feynman-Kac formula
\begin{equation}\label{eq:CHNSBSDE}
	\left\{
	\begin{aligned}
		& dX_{s}^{u}=\sqrt{2\nu}dW_{s}^{u}, \quad  s \in [0,T],\\
		& dX_{s}^{\hat{\phi}}=\sqrt{2\lambda_1}dW_{s}^{\hat{\phi}}, \quad  s \in [0,T],\\
		& dX_{s}^{\hat{\mu}}=\sqrt{2\lambda_2}dW_{s}^{\hat{\mu}}, \quad  s \in [0,T],\\
		& (X_{0}^{u},X_{0}^{\hat{\phi}}, X_{0}^{\hat{\mu}})^T = (x_0, x_0, x_0)^T,\\
		& -dY_{s}^{u}=F_{s}^{u}ds-\sqrt{2\nu}(Z_{s}^{u})^TdW_{s}^{u} , \quad  s \in [0,T],\\
		& -dY_{s}^{\hat{\phi}}=\hat{F}_{s}^{\hat{\phi}}ds-\sqrt{2\lambda_1}(Z_{s}^{\hat{\phi}})^TdW_{s}^{\hat{\phi}} , \quad  s \in [0,T],\\
		&  -dY_{s}^{\hat{\mu}}=\hat{F}_{s}^{\hat{\mu}}ds-\sqrt{2\lambda_2}(Z_{s}^{\hat{\mu}})^TdW_{s}^{\hat{\mu}}, \quad  s \in [0,T],\\
		& u(T) = g_{u}, \quad \phi(T) = g_{\phi},
	\end{aligned}
	\right.
\end{equation}
with $ F^u = (u\cdot\nabla)u+\nabla p+C\phi\nabla\mu+f_1$ and $$ (\hat{F}^{\hat{\phi}}, \hat{F}^{\hat{\mu}})^T = R^{-1}\left(f_2 + u\cdot\nabla\phi, \frac{S}{L_d}\left(u\cdot\nabla\phi + f_2\right)-\frac{1}{\delta}\left(\mu + \phi-\phi^3\right)\right)^T. $$ The Euler scheme of \eqref{eq:CHNSBSDE} for $n=0,\cdots, N-1$ is 
\begin{equation}\label{eq:CHNSEuler1}
\left\{
	\begin{aligned}
		& X_{t_{n+1}}^u=X_{t_{n}}^u+\sqrt{2\nu}\Delta W_{t_{n}}^u,\\
		& X_{t_{n+1}}^{\hat{\phi}}=X_{t_{n}}^{\hat{\phi}}+\sqrt{2\lambda_1}\Delta W_{t_{n}}^{\hat{\phi}},\\
		& X_{t_{n+1}}^{\hat{\mu}}=X_{t_{n}}^{\hat{\mu}}+\sqrt{2\lambda_2}\Delta W_{t_{n}}^{\hat{\mu}},\\
		&\widetilde{Y}_{t_{n+1}}^u(X_{t_{n+1}}^u)=Y_{t_{n}}^u(X_{t_{n}}^u)+\sqrt{2\nu}(Z_{t_{n}}^u(X_{t_{n}}^u))^T\Delta W_{t_{n}}^u\\
            &-F_{t_{n}}^u(X_{t_{n}}^u,\nabla P_{t_{n}}(X_{t_{n}}^u),Y_{t_{n}}^u(X_{t_{n}}^u),Y_{t_{n}}^{\phi}(X_{t_{n}}^u),Z_{t_{n}}^{\mu}(X_{t_{n}}^u))\Delta t_{n},\\
		&\widetilde{Y}_{t_{n+1}}^{\hat{\phi}}(X_{t_{n+1}}^{\hat{\phi}})=Y_{t_{n}}^{\hat{\phi}}(X_{t_{n}}^{\hat{\phi}})+\sqrt{2\lambda_1}(Z_{t_{n}}^{\hat{\phi}}(X_{t_{n}}^{\hat{\phi}}))^T\Delta W_{t_{n}}^{\hat{\phi}}\\
        &-\hat{F}_{t_{n}}^{\hat{\phi}}(X_{t_{n}}^{\hat{\phi}},Y_{t_{n}}^u(X_{t_{n}}^{\hat{\phi}}),Y_{t_{n}}^{\phi}(X_{t_{n}}^{\hat{\phi}}),Y_{t_{n}}^{\mu}(X_{t_{n}}^{\hat{\phi}}),Z_{t_{n}}^{\phi}(X_{t_{n}}^{\hat{\phi}}))\Delta t_{n},\\
	&\widetilde{Y}_{t_{n+1}}^{\hat{\mu}}           
        (X_{t_{n+1}}^{\hat{\mu}})=Y_{t_{n}}^{\hat{\mu}}(X_{t_{n}}^{\hat{\mu}})+\sqrt{2\lambda_2}(Z_{t_{n}}^{\hat{\mu}}(X_{t_{n}}^{\hat{\mu}}))^T\Delta W_{t_{n}}^{\hat{\mu}}\\
        &-\hat{F}_{t_{n}}^{\hat{\mu}}(X_{t_{n}}^{\hat{\mu}},Y_{t_{n}}^u(X_{t_{n}}^{\hat{\mu}}),Y_{t_{n}}^{\phi}(X_{t_{n}}^{\hat{\mu}}),Y_{t_{n}}^{\mu}(X_{t_{n}}^{\hat{\mu}}),Z_{t_{n}}^{\phi}(X_{t_{n}}^{\hat{\mu}}))
        \Delta t_{n},
	\end{aligned}
 \right.
\end{equation}
where  $\Delta t_{n}=t_{n+1}-t_{n}=\frac{T}{N} =\delta ,\Delta W_{t_{n}}^{u} = W_{t_{n+1}}^{u}-W_{t_{n}}^{u}, \Delta W_{t_{n}}^{\hat{\phi}} = W_{t_{n+1}}^{\hat{\phi}}-W_{t_{n}}^{\hat{\phi}}$ and $\Delta W_{t_{n}}^{\hat{\mu}} = W_{t_{n+1}}^{\hat{\mu}}-W_{t_{n}}^{\hat{\mu}}$. 
The $Y_{t_{n}}^u(X_{t_{n}}^u)$, $Y_{t_{n}}^{\hat{\phi}}(X_{t_{n}}^{\hat{\phi}})$ and $Y_{t_{n}}^{\hat{\mu}}(X_{t_{n}}^{\hat{\mu}})$ represent the estimated values of $u(t_n,X_{t_{n}}^u)$, $\hat{\phi}(t_n,X_{t_{n}}^{\hat{\phi}})$ and $\hat{\mu}(t_n,X_{t_{n}}^{\hat{\mu}})$, respectively.
The $\widetilde{Y}_{t_{n+1}}^u(X_{t_{n+1}}^u)$, $\widetilde{Y}_{t_{n+1}}^{\hat{\phi}}(X_{t_{n+1}}^{\hat{\phi}})$ and $\widetilde{Y}_{t_{n+1}}^{\hat{\mu}}(X_{t_{n+1}}^{\hat{\mu}})$ are the reference values of $Y_{t_{n+1}}^u(X_{t_{n+1}}^u)$, $Y_{t_{n+1}}^{\hat{\phi}}(X_{t_{n+1}}^{\hat{\phi}})$ and $Y_{t_{n+1}}^{\hat{\mu}}(X_{t_{n+1}}^{\hat{\mu}})$, respectively, which are obtained from the last three equations in \eqref{eq:CHNSEuler1}.

The $(Y_{t_{n}}^u,P_{t_{n}})^T$ represents the estimated value of $(u,p)^T$ at time $t_n$ given by the FNN $\mathcal{U}_{\theta_1}$. The $(Y_{t_{n}}^{\phi}, Y_{t_{n}}^{\mu})^T$ represents the estimated value of $(\phi,\mu)^T$ at time $t_n$ given by the FNN $\mathcal{U}_{\theta_2}$.    
The calculations of $Y_{t_{n}}^{\hat{\phi}}(X_{t_{n}}^{\hat{\phi}})$ and $Y_{t_{n}}^{\hat{\mu}}(X_{t_{n}}^{\hat{\mu}})$ are based on \eqref{eq:diag}. We choose $K$ different initial sampling points for training.
The algorithm of the proposed scheme is summarized as Algorithm \ref{alg:Algorithm 3}.
\begin{algorithm}
	\caption{ Algorithm  for  the
		Cahn-Hilliard-Navier-Stokes equation}
	\label{alg:Algorithm 3}
	\begin{algorithmic}[1]
		\Require
		 Number of initial sampling points $K$, terminal time $T$, number of time intervals $N$, viscosity coefficient $\nu$,  diffusion coefficient $L_d$, parameters $\gamma$, $\delta$, $C$ and $S$, matrix $R^{-1}$, eigenvalues $\lambda_{1,2}$,  maximum number of training iterations $M_{iter}$.
		\Ensure
		The optimal FNNs $\mathcal{U}_{\theta_1}$ and $\mathcal{U}_{\theta_2}$.
		\State Initialize the FNNs $\mathcal{U}_{\theta_1}$ and $\mathcal{U}_{\theta_2}$;
        \State  Select initial sampling points $x_0$ by uniform distribution;
		\State Generate independent $d$-dimensional standard Brownian motions $W_{t_n}^u
		$, $W_{t_n}^{\hat{\phi}}$ and $W_{t_n}^{\hat{\mu}}(n = 0,...,N)$;
		\State  Compute $X_{t_{n+1}}^u$, $X_{t_{n+1}}^{\hat{\phi}}$ and $X_{t_{n+1}}^{\hat{\mu}}$ according to \eqref{eq:CHNSEuler1} for $n = 0,..., N-1$;
        \State According to \eqref{eq:diag} and \eqref{eq:CHNSEuler1}, use the FNNs  $\mathcal{U}_{\theta_1}$ and $\mathcal{U}_{\theta_2}$ with AD to calculate
        $\widetilde{Y}_{t_{n+1}}^u(X_{t_{n+1}}^{u})$, $\widetilde{Y}_{t_{n+1}}^{\hat{\phi}}(X_{t_{n+1}}^{\hat{\phi}})$ and $\widetilde{Y}_{t_{n+1}}^{\hat{\mu}}(X_{t_{n+1}}^{\hat{\mu}})$;
		\State Minimize the loss function by the Adam algorithm
        \begin{equation}\label{eq:lossCHNSBSDE1}
    	\begin{aligned}
		&l(\theta_1,\theta_2) = \frac{1}{K}\sum_{k=1}^{K}\left[\frac{1}{N}\sum_{n=0}^{N-1}\left(\vert (\widetilde{Y}_{t_{n+1}}^u-Y_{t_{n+1}}^u)(X_{t_{n+1},k}^{u})\vert^{2}\right.\right.\\
        &\left.+\vert (\widetilde{Y}_{t_{n+1}}^{\hat{\phi}}-Y_{t_{n+1}}^{\hat{\phi}})(X_{t_{n+1},k}^{\hat{\phi}})\vert^{2}+\vert (\widetilde{Y}_{t_{n+1}}^{\hat{\mu}}-Y_{t_{n+1}}^{\hat{\mu}})(X_{t_{n+1},k}^{\hat{\mu}})\vert^{2}\right)
        \\
        &+ \alpha_1\vert (g_u-Y_{t_N}^u)(X_{t_{N},k}^{u\cup\hat{\phi}\cup\hat{\mu}})\vert^{2}+\alpha_2\vert (g_{\phi}-Y_{t_N}^{\phi})(X_{t_{N},k}^{u\cup\hat{\phi}\cup\hat{\mu}})\vert^{2}
        \\
        &\left.+\frac{\alpha_3}{N+1}\sum_{n=0}^{N}\vert\nabla\cdot Y_{t_{n}}^u(X_{t_{n},k}^{u\cup\hat{\phi}\cup\hat{\mu}})\vert^2\right],\\
        \end{aligned}
    \end{equation}
    where $X_{t_{n},k}^{u\cup\hat{\phi}\cup\hat{\mu}} = X_{t_n,k}^u\cup X_{t_{n},k}^{\hat{\phi}}\cup X_{t_{n},k}^{\hat{\mu}}$ and $\alpha_i,i=1,2,3$ are the weights of the components of the loss function. The subscript $k$ corresponds to the $k$-th initial sampling point;
		\State Repeat procedures 2 to 6 until the maximum number of training iterations $M_{iter}$ is reached.
	\end{algorithmic}
\end{algorithm}
For the CHNS equation with the Dirichlet, Neumann and periodic boundary conditions, the similar methods as in Section \ref{sec03} and Section \ref{sec04} can be used.

\section{Numerical experiments}
\label{sec06}
In this section, we present a series of numerical results to validate our methods. For quantitative comparison, we calculate the error of the numerical solution $Y_{t_0}$ and the exact solution $u_{t_0}$ in the relative $L^{\infty}$ norm and relative $L^2$ norm, which are defined as  
$$\vert\vert e\vert\vert_{L^{\infty}} = \frac{\vert\vert Y_{t_0}-u_{t_0}\vert\vert_{L^{\infty}}}{\vert\vert u_{t_0}\vert\vert_{L^{\infty}}}, \ \ \ \vert\vert e\vert\vert_{L^2}=\frac{\vert\vert Y_{t_0}-u_{t_0}\vert\vert_{L^{2}}}{\vert\vert u_{t_0}\vert\vert_{L^{2}}}.$$ The total number of training iterations is given by 1E+5, which is divided into 2E+4, 3E+4, 3E+4 and 2E+4 iterations with learning rates of 5E-3, 5E-4, 5E-5 and 5E-6, respectively, as the way in \cite{raissi2024forward}. We employ the Adam
optimizer to train FNNs. For each training step, we train the FNNs using 100 points randomly selected by the Latin hypercube sampling technique (LHS) in the domain. After the training process, we randomly pick 10000 points by the same method in the domain to test the FNNs. We set 4 hidden layers for  the FNNs and each hidden layer consists of 30 neurons. The cosine function is taken as the activation function of the FNNs if we do not specify otherwise. In each numerical example, we use a set of the time interval $\Delta t$, weights $\alpha_i$ and the stabilization parameter $S$. How to select or adjust these hyperparameters will be our future work. In our simulations, we use AMD Ryzen 7 3700X CPU and NVIDIA GTX 1660 SUPER GPU to train FNNs. The parameters and settings of numerical experiments are summarized in the Table \ref{Table:set}.
\begin{table}[htp]
	\begin{center} 
\begin{tabular}{ccc}
\hline 
\text { total number of training iterations } & 1E+5 \\
\text { number of iterations per segment } & [2E+4, 3E+4, 3E+4, 2E+4]   \\
\text { learning rate per segment } & [5E-3, 5E-4, 5E-5, 5E-6]   \\
\text { optimization algorithm }  & Adam \\
\text { structure of neural networks} & [30,30,30,30]  \\
\text { number of training points } & 100 \\
\text { number of test points } & 10000 \\
\text { point selection method }  & \text { LHS } \\
\text { activation function } &  cos \\
\text { CPU } & \text {AMD Ryzen 7 3700X} \\
\text { GPU } & \text {NVIDIA GTX 1660 SUPER} \\
\hline
\end{tabular}
\caption{
Parameters and settings for numerical experiments }\label{Table:set}
	\end{center}
\end{table}
\subsection{Navier-Stokes equation}
\label{sec06:sec01}
\begin{figure}[htp]
	\centering 
	\subfigure[Relative $L^2$ error]{		 
		\includegraphics[scale=0.45]{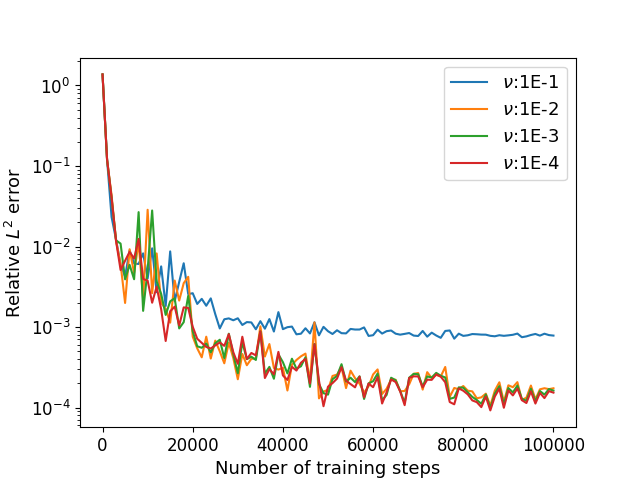}}
	\subfigure[Training loss]{		 
		\includegraphics[scale=0.45]{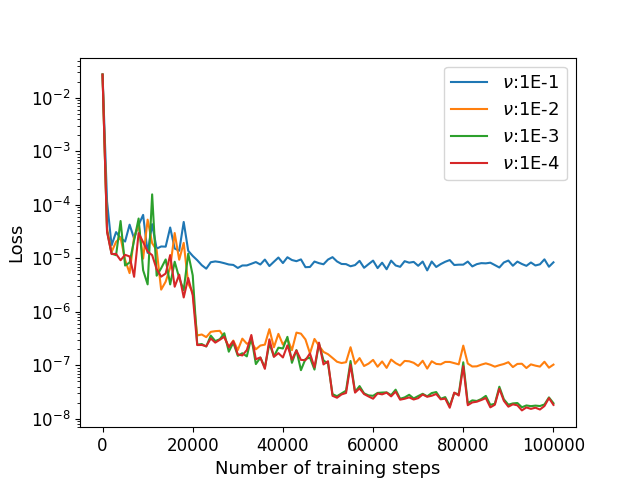}}
	\caption{Relative $L^{2}$ errors of $\mathbf{u}$ and training losses of the Algorithm \ref{alg:Algorithm 1} for \eqref{eq:Taysolu1} with  different $\nu$.}
	\label{fig:fig3}
\end{figure}

In this section, we numerically simulate the Taylor-Green vortex flow, which is a classical model to test numerical schemes for the incompressible Navier-Stokes equation. First, we consider the explicit 2D Taylor-Green vortex flow
\begin{equation}\label{eq:Taysolu1}
	\left\{
	\begin{aligned}
		&u_1(t,x)=-\cos(x_1)\sin(x_2)\exp(-2\nu t),\\
		&u_2(t,x)=\sin(x_1)\cos(x_2)\exp(-2\nu t),\\
		&p(t,x)=-\frac{1}{4}\left(\cos(2x_1)+\cos(2x_2)\right)\exp(-4\nu t)+c,
	\end{aligned}
	\right.
\end{equation}
for $(t,x)= (0,T]\times[0,2\pi]^2$ with constant $c\in \mathbb{R}$ and initial condition
\begin{equation}\label{eq:Taysolu1init}
	\left\{
	\begin{aligned}
		&u_1(0,x)=-\cos(x_1)\sin(x_2),\\
		&u_2(0,x)=\sin(x_1)\cos(x_2).
	\end{aligned}
	\right.
\end{equation}
Algorithm \ref{alg:Algorithm 1} is employed to estimate $\mathbf{u}(T,x) = (u_1,u_2)^T$  with $T = 0.1$, $f = 0$, $ N=5 $, $\Delta t = \frac{T}{N}=0.02$ and $\alpha_1=\alpha_2=0.1$.  The numerical results of the  errors for $\mathbf{u}$ and $\nabla p$ with different viscosity $\nu$ are shown in Table \ref{Table:1}.
The relative $L^{2}$ errors and the training losses with different training steps are shown in Figure \ref{fig:fig3}. It is observed that these values decrease with parameter $\nu$ decreases. Similar phenomena will occur in the later experiments. Our method is not sensitive to the parameter $\nu$.
The training time is 500s for each case, which is a acceptable cost. 

\begin{table}[htp]
	\begin{center} 
		\begin{tabular}{|l|l|l|l|l|l|l|l|l|}
			\hline
			$\nu$ & 1E-1 & 1E-2  & 1E-3 & 1E-4         \\
			\hline  
			$\vert\vert e_1\vert\vert_{L^{\infty}} $  & 1.60E-2 & 8.58E-3 & 7.41E-3 & 6.92E-3    \\
			\hline 
			$\vert\vert e_2\vert\vert_{L^{\infty}} $  & 1.66E-2 & 5.72E-3 & 6.30E-3 & 6.45E-3   \\
			\hline
			$\vert\vert e_1\vert\vert_{L^2} $  & 7.23E-3 & 2.00E-3 & 9.02E-4 & 8.52E-4   \\
			\hline
			$\vert\vert e_2\vert\vert_{L^2} $  & 8.61E-3 & 1.39E-3 & 8.59E-4 & 8.36E-4   \\
			\hline
                $\vert\vert e_{\nabla p}\vert\vert_{L^2} $  & 1.12E-1 & 2.41E-2 & 9.45E-3 & 8.40E-3   \\
                \hline
                time  & \multicolumn{4}{|c|}{500s}    \\
			\hline
		\end{tabular}
		\caption{Relative $L^{\infty}$ and  $L^2$ errors for \eqref{eq:Taysolu1} performed by the Algorithm \ref{alg:Algorithm 1} with $\nu=$ 1E-1,\ 1E-2,\ 1E-3 and 1E-4. }\label{Table:1}
	\end{center}
\end{table}

\begin{figure}[htp]
	\centering 
	\subfigure[Relative $L^2$ error]{		 
		\includegraphics[scale=0.45]{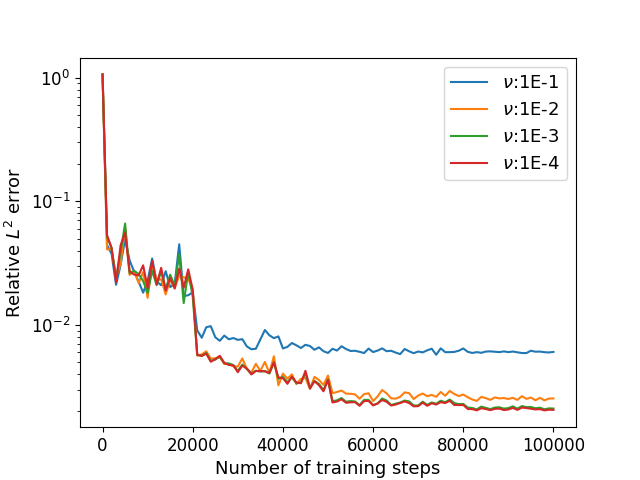}}
	\subfigure[Training loss]{		 
		\includegraphics[scale=0.45]{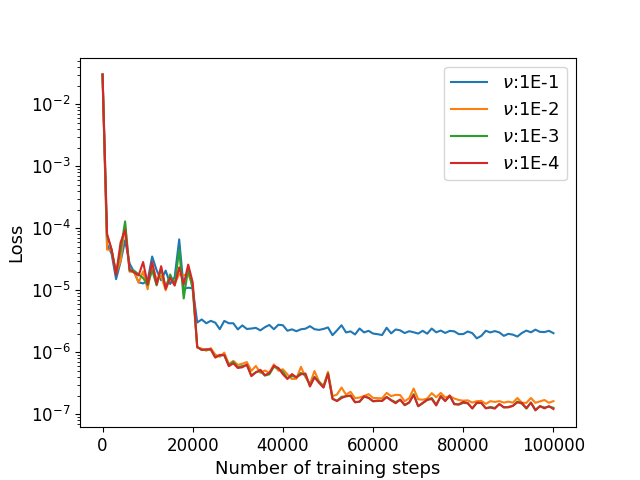}}
	\caption{Relative $L^{2}$ errors and training losses of the Algorithm \ref{alg:Algorithm 1} for \eqref{eq:ABCflow} with  different $\nu$.} 
	\label{fig:fig4}
\end{figure}

The 3D Arnold-Beltrami-Childress (ABC) flow is as follows
\begin{equation}\label{eq:ABCflow}
	\left\{
	\begin{aligned}
		& u_1(t,x)=(A\sin(x_3)+C\cos(x_2))e^{-\nu t},\\
		& u_2(t,x)=(B\sin(x_1)+A\cos(x_3))e^{-\nu t},\\
		& u_3(t,x)=(C\sin(x_2)+B\cos(x_1))e^{-\nu t},\\
		& p(t,x)=-(BC\cos(x_1)\sin(x_2)+AB\sin(x_1)\cos(x_3)\\
		&\quad\quad\quad\quad+AC\sin(x_3)\cos(x_2))e^{-2\nu t}+c,
	\end{aligned}
	\right.
\end{equation}
for $(t,x)=(0,T]\times[0,2\pi]^3$ with parameters $A,B,C\in\mathbb{R}$, constant $c\in\mathbb{R}$ and initial condition
\begin{equation}\label{eq:ABCflowinit}
	\left\{
	\begin{aligned}
		& u_1(0,x)=(A\sin(x_3)+C\cos(x_2)),\\
		& u_2(0,x)=(B\sin(x_1)+A\cos(x_3)),\\
		& u_3(0,x)=(C\sin(x_2)+B\cos(x_1)).
	\end{aligned}
	\right.
\end{equation}
We estimate $\mathbf{u}(T,x) = (u_1,u_2,u_3)^T$ by applying the Algorithm \ref{alg:Algorithm 1} with parameters $A=B=C=0.5$, $T = 0.1$, $f = 0$,  $ N=5 $, $\Delta t = \frac{T}{N}=0.02$, $\alpha_1=\alpha_2=0.1$. The numerical results of the  errors for $\mathbf{u}$ and $\nabla p$ with different viscosity $\nu$ are shown in Table \ref{Table:2}. The relative $L^{2}$ errors and the training losses with different training steps are shown in Figure \ref{fig:fig4}. The training time is 700s for each case, which is not too much longer than the 2D simulations. 
\begin{table}[htp]
	\begin{center} 
		\begin{tabular}{|l|l|l|l|l|l|l|l|l|}
			\hline
			$\nu$ & 1E-1 & 1E-2  & 1E-3 & 1E-4         \\
			\hline
			$\vert\vert e_1\vert\vert_{L^{\infty}} $   & 1.46E-2 & 8.77E-3 & 9.19E-3 & 9.16E-3   \\
			\hline 
			$\vert\vert e_2\vert\vert_{L^{\infty}} $   & 1.41E-2 & 8.12E-3 & 9.34E-3 & 9.60E-3   \\
			\hline
			$\vert\vert e_3\vert\vert_{L^{\infty}} $   & 1.24E-2 & 9.32E-3 & 8.83E-3 & 8.81E-3   \\
			\hline
			$\vert\vert e_1\vert\vert_{L^2} $  & 6.20E-3 & 2.01E-3 & 1.81E-3 & 1.83E-3   \\
			\hline
			$\vert\vert e_2\vert\vert_{L^2} $  & 4.60E-3 & 2.36E-3 & 2.08E-3 & 2.05E-3   \\
			\hline
			$\vert\vert e_3\vert\vert_{L^2} $  & 7.05E-3 & 3.12E-3 & 2.38E-3 & 2.26E-3   \\
			\hline
                $\vert\vert e_{\nabla p}\vert\vert_{L^2} $  & 1.70E-1 & 6.92E-2 & 5.72E-2 & 5.58E-2   \\
			\hline
                time  & \multicolumn{4}{|c|}{700s}    \\
			\hline
		\end{tabular}
		\caption{Relative $L^{\infty}$ and $L^2$ errors for \eqref{eq:ABCflow} performed by the Algorithm \ref{alg:Algorithm 1} with $\nu=$1E-1,\ 1E-2,\ 1E-3 and 1E-4. }\label{Table:2}
	\end{center}
\end{table}

\begin{figure}[htp]
	\centering 
	\subfigure[Relative $L^2$ error]{		 
		\includegraphics[scale=0.45]{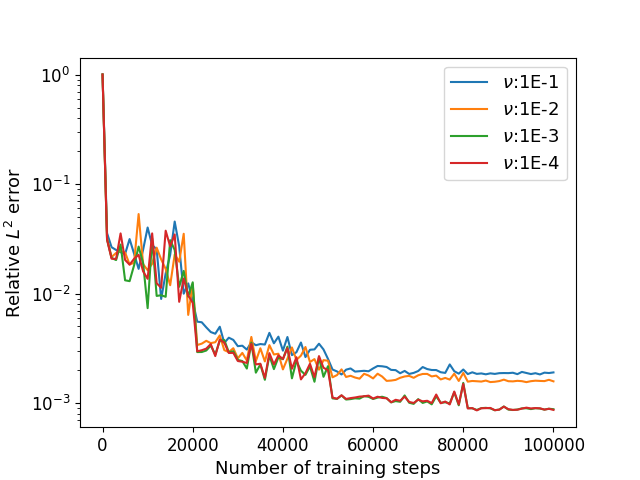}}
	\subfigure[Training loss]{		 
		\includegraphics[scale=0.45]{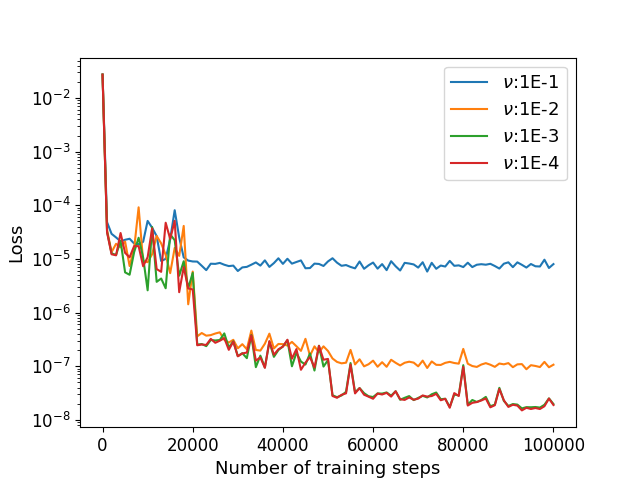}}
	\caption{Relative $L^{2}$ errors and training losses of the algorithm in Section \ref{sec03:sec03} for \eqref{eq:Taysolu1} with the Dirichlet boundary condition and different $\nu$.}
	\label{fig:fig5}
\end{figure}

\begin{figure}[htp]
	\centering 
	\subfigure[Relative $L^2$ error]{		 
		\includegraphics[scale=0.45]{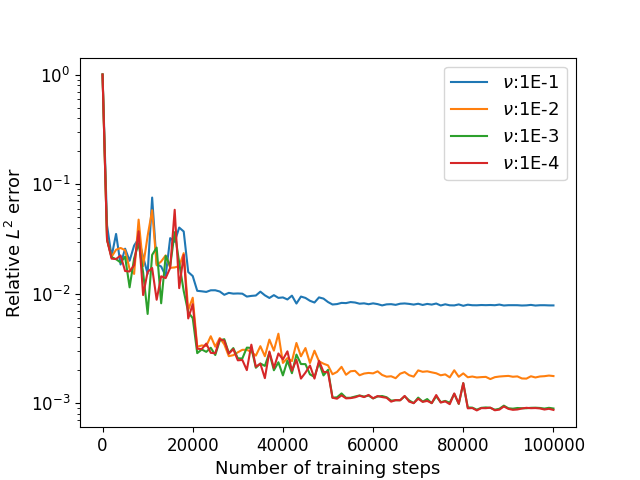}}
	\subfigure[Training loss]{		 
		\includegraphics[scale=0.45]{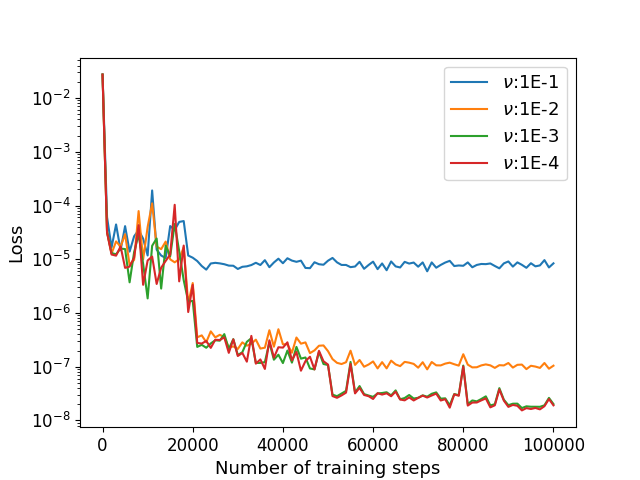}}
	\caption{Relative $L^{2}$ errors and training losses of the algorithm in Section \ref{sec03:sec04} for \eqref{eq:Taysolu1} with the Neumann boundary condition and different $\nu$.}
	\label{fig:fig6}
\end{figure}

 Next, we consider the 2D Taylor-Green vortex flow \eqref{eq:Taysolu1}--\eqref{eq:Taysolu1init} with the Dirichlet boundary condition \eqref{eq:dboundary} for $(t,x)=(0,T]\times(0,2\pi)^2$.
 The other parameters remain the same as the first example. We use the algorithm in Section \ref{sec03:sec03} and the numerical results of the errors for $\mathbf{u}$ and $\nabla p$ with different viscosity parameters $\nu$ are shown in Table \ref{Table:3} and Figure \ref{fig:fig5} depicts the training processes. We also consider the 2D Taylor-Green vortex flow \eqref{eq:Taysolu1}--\eqref{eq:Taysolu1init} with the Neumann boundary condition \eqref{eq:nboundary} for $(t,x)=(0,T]\times(0,2\pi)^2$.
 We use the algorithm in Section \ref{sec03:sec04} and the numerical results of the errors for $\mathbf{u}$
and $\nabla p$ with different viscosity parameters $\nu$ are shown in Table \ref{Table:4} and the training processes are shown in Figure \ref{fig:fig6}. 
\begin{table}[htp]
	\begin{center} 
		\begin{tabular}{|l|l|l|l|l|l|l|l|l|}
			\hline
			$\nu$ & 1E-1 & 1E-2  & 1E-3 & 1E-4         \\
			\hline 
			$\vert\vert e_1\vert\vert_{L^{\infty}} $  & 4.96E-3 & 8.19E-3 & 7.13E-3 & 6.99E-3   \\
			\hline 
			$\vert\vert e_2\vert\vert_{L^{\infty}} $  & 5.44E-3 & 5.70E-3 & 6.09E-3 & 6.41E-3  \\
			\hline  
			$\vert\vert e_1\vert\vert_{L^2}$ & 1.89E-3 & 1.57E-3 & 8.90E-4 & 8.70E-4   \\
			\hline
			$\vert\vert e_2\vert\vert_{L^2} $ & 1.92E-3 & 1.59E-3 & 8.62E-4 & 8.57E-4   \\
			\hline
                $\vert\vert e_{\nabla p}\vert\vert_{L^2} $  & 3.26E-2 & 2.20E-2 & 9.28E-3 & 8.83E-3   \\
			\hline
                time  & \multicolumn{4}{|c|}{700s}\\
                \hline
		\end{tabular}
		\caption{Relative $L^{\infty}$ and $L^2$ errors  for \eqref{eq:Taysolu1} with the Dirichlet boundary condition performed by the algorithm in Section \ref{sec03:sec03} with $\nu=$ 1E-1,\ 1E-2,\ 1E-3 and 1E-4. }\label{Table:3}
	\end{center}
\end{table}

\begin{table}[htp]
	\begin{center} 
		\begin{tabular}{|l|l|l|l|l|l|l|l|l|}
			\hline
			$\nu$ & 1E-1 & 1E-2  & 1E-3 & 1E-4         \\
                \hline 
			$\vert\vert e_1\vert\vert_{L^{\infty}} $  & 1.50E-2 & 9.45E-3 & 7.65E-3 & 7.35E-3    \\
			\hline 
			$\vert\vert e_2\vert\vert_{L^{\infty}} $  & 1.55E-2 & 6.10E-3 & 6.50E-3 & 6.64E-3   \\
			\hline
			$\vert\vert e_1\vert\vert_{L^2}$  & 8.05E-3 & 1.98E-3 & 9.27E-4 & 8.83E-4   \\
			\hline
			$\vert\vert e_2\vert\vert_{L^2} $  & 7.59E-3 & 1.52E-3 & 8.60E-4 & 8.53E-4  \\
			\hline
                $\vert\vert e_{\nabla p}\vert\vert_{L^2} $  & 1.11E-1 & 2.43E-2 & 9.36E-3 & 9.00E-3   \\
			\hline
                time  & \multicolumn{4}{|c|}{900s} \\
                \hline
		\end{tabular}
		\caption{Relative $L^{\infty}$ and $L^2$ errors  for \eqref{eq:Taysolu1} with the Neumann boundary  condition performed by the algorithm in Section \ref{sec03:sec04} with $\nu=$ 1E-1,\ 1E-2,\ 1E-3 and 1E-4. }\label{Table:4}
	\end{center}
\end{table}

\begin{figure}[htp]
	\centering 
	\subfigure[stream function]{		 
		\includegraphics[scale=0.45]{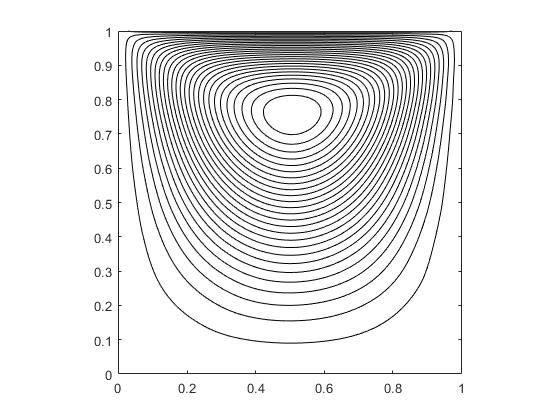}}
	\subfigure[pressure $p$]{		 
		\includegraphics[scale=0.45]{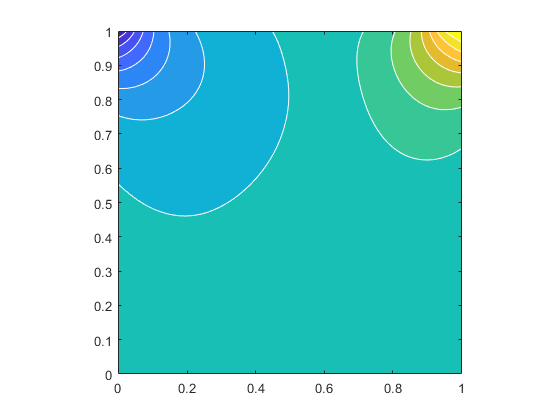}}
	\caption{The stream function and pressure $p$  for 2D lid driven cavity flow at $T = 0.5$ with $\nu = 0.1$.}
	\label{fig:fig13}
\end{figure}

Now, we consider the 2D lid driven cavity flow for $(t,x)=(0,T]\times(0,1)^2$ with the
boundary and initial conditions
\begin{equation}\label{eq:cavity}
	\left\{
	\begin{aligned}
		& \mathbf{u}(t,x)=(1,0), \quad  x\in\partial \Omega_{u}\\
		& \mathbf{u}(t,x)=(0,0), \quad  x\in\partial\Omega\backslash\partial \Omega_{u}\\
		& \mathbf{u}(0,x)=(0,0), \quad x \in \Omega,
	\end{aligned}
	\right.
\end{equation}
where $\partial \Omega_{u}$ represents the upper boundary. We  utilize the algorithm in Section \ref{sec03:sec03} to  simulate $\mathbf{u}$. 
We impose boundary conditions $u_i=0, i=1,2$ to the network in the training process and let \begin{equation}
      \left\{
    \begin{aligned}\label{eq:enforce1}
     &Y_{t_{n}}^1(x_1,x_2) = 8x_1(x_1-1)x_2Y_{t_{n}}^1(x_1,x_2),\\
     &Y_{t_{n}}^2(x_1,x_2) = 8x_1(x_1-1)x_2(x_2-1)Y_{t_{n}}^2(x_1,x_2),
    \end{aligned} \right.
     \end{equation} 
where $Y_{t_{n}}^i(x_1,x_2)$  represents the estimate of $u_i(t_n,x_1,x_2)$ output by the FNN. Therefore, it is easily verified that $Y_{t_{n}}$ satisfies the boundary conditions of $u_i=0, i=1,2$. For the condition of $u_1=1$ on $\partial \Omega_{u}$ , we add the following additional term to the loss function in Algorithm \ref{alg:Algorithm 1}
\begin{equation}\label{eq:leftandright1}
	\begin{aligned} 
		& \frac{\alpha_3}{K_u}\sum_{k=1}^{K_u}\left[\frac{1}{N+1}\sum_{n=0}^{N}\vert Y_{t_{n}}^1(X^u_{t_n,k})+1\vert^2\right],
	\end{aligned}
\end{equation}
where $X^u_{t_n,k}$ represents the $k$-th point among the $K_u$ points selected on $\partial \Omega_{u}$ at time $t_n$. The parameters are chosen as $T = 0.5$, $f = 0$, $ \nu = 0.1 $, $ N = 25 $, $K_u = 100$, $\Delta t = 0.02$,  $\alpha_1=\alpha_2=\alpha_3=0.01$. To improve accuracy and save training time, we  use the time adaptive approach \RNum{2} mentioned in \cite{wight2020solving}.
At $T = 0.5$, the stream function and the pressure $p$ with $\nu = 0.1$ are visually shown in Figure \ref{fig:fig13}. These results are consistent with benchmark results. 

\begin{figure}[htp]
	\centering 
	\subfigure{		 
		\includegraphics[scale=0.9]{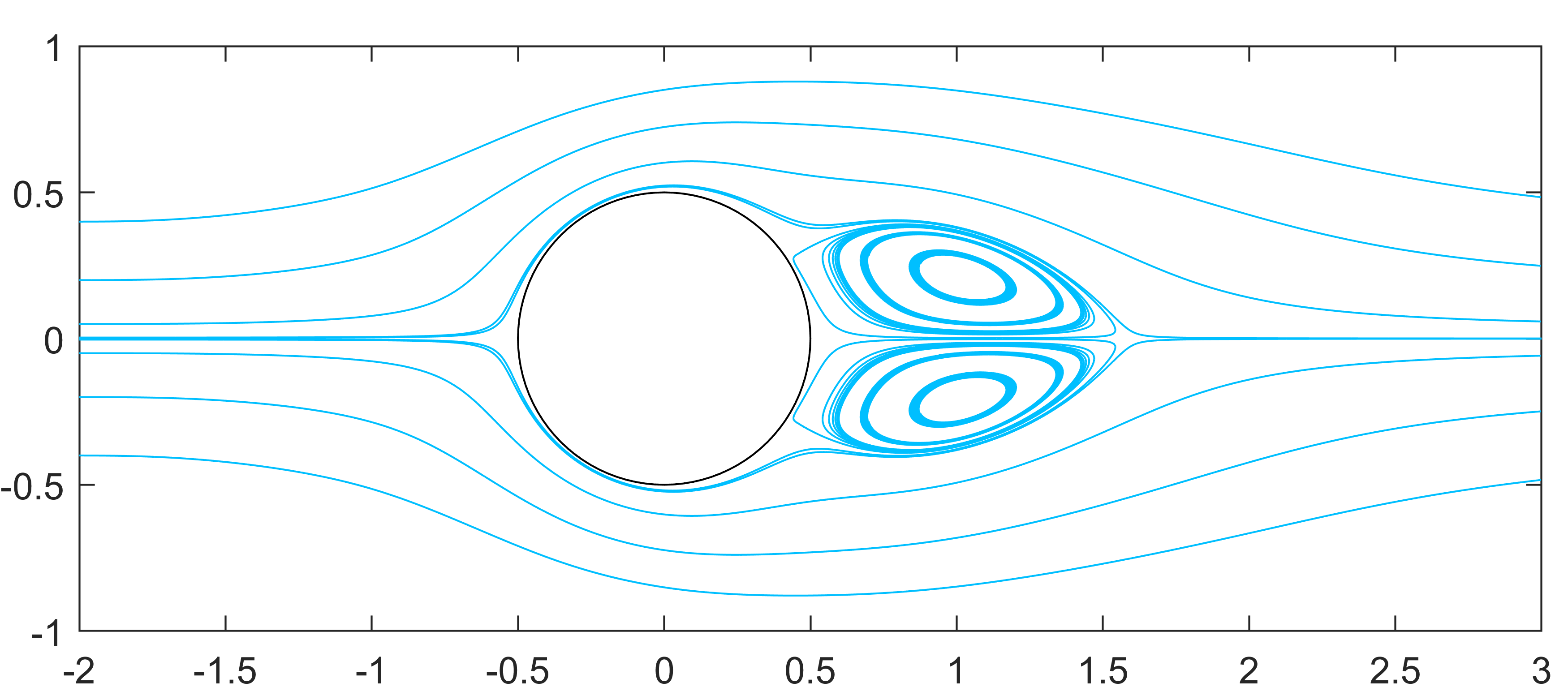}}
	\caption{The streamline for the flow past a circular obstacle at $T = 1.0$ with $\nu = 0.025$.}
	\label{fig:fig14}
\end{figure}

Finally, we consider that the flow past a circular obstacle for $(t,x)=(0,T]\times(-2,10)\times(-2,2)$. The center of the obstacle is at position $(0,0)$ with the diameter $D=1$. The boundary and initial conditions are given 
\begin{equation}\label{eq:fpco}
	\left\{
	\begin{aligned}
		& u_2(t,x) = 0, \quad  x\in\partial \Omega_{u}\cup\partial \Omega_{d},\\\
		& \mathbf{u}(t,x)=(0,0), \quad  x\in\partial \Omega_{c},\\
            & \mathbf{u}(t,x)=(u_{in},0), \quad  x\in\partial \Omega_{l},\\
            & p\mathbf{n}-\nu\nabla \mathbf{u}\cdot \mathbf{n}=0, \quad  x\in\partial \Omega_{r},\\
		& \mathbf{u}(0,x)=(u_{in},0), \quad x\in\Omega,
	\end{aligned}
	\right.
\end{equation}
where $\partial \Omega_{u},\partial \Omega_{d},\partial\Omega_{l},\partial \Omega_{r},\partial \Omega_{c}$ represent the upper, lower, left, right boundaries and the surface of the obstacle. The $u_{in}$ is the inlet velocity and $\mathbf{n}$ is the unit normal vector at $\partial\Omega$ pointing outward of $\Omega$.  We  utilize the algorithm in Section \ref{sec03:sec03} to  deal with the Dirichlet boundary conditions on $\partial\Omega_{u}\cup\partial\Omega_{d}\cup\partial\Omega_{l}\cup\partial\Omega_{c}$, while we  utilize the algorithm in Section \ref{sec03:sec04} to  deal with the condition on $\partial\Omega_{r}$. 
We let \begin{equation}\label{eq:enforce2}
\left\{
    \begin{aligned}
    &Y_{t_{n}}^1(x_1,x_2) = \frac{8Y_{t_{n}}^1(x_1,x_2)(x_1^2+x_2^2-0.25)}{(x_1^2+x_2^2)},\\
    &Y_{t_{n}}^2(x_1,x_2) = \frac{8Y_{t_{n}}^2(x_1,x_2)(x_1^2+x_2^2-0.25)(x_1+2)(x_2-2)(x_2+2)}{(x_1^2+x_2^2)},\\
    \end{aligned}
    \right.
\end{equation} 
where $Y_{t_{n}}^i(x_1,x_2)$  represents the estimate of $u_i(t_n,x_1,x_2)$ output by the FNN. Therefore, it is easily verified that $Y_{t_{n}}$ satisfies the boundary conditions of $u_i=0, i=1,2$.
For the conditions of $u_1=u_{in}$ on $\partial\Omega_{l}$ and $p\mathbf{n}-\nu\nabla \mathbf{u}\cdot \mathbf{n}=0$ on $\partial\Omega_{r}$, we add the following additional term to the loss function in Algorithm \ref{alg:Algorithm 1}
\begin{equation}\label{eq:leftandright2}
	\begin{aligned} 
		& \frac{\alpha_3}{N+1}\sum_{n=0}^{N}\left[\frac{1}{K_l}\sum_{k=1}^{K_l}\vert Y_{t_{n}}^1(X^l_{t_n,k})+u_{\infty}\vert^2+\frac{1}{K_r}\sum_{k=1}^{K_r}\vert (P_{t_{n}}\mathbf{n}+\nu Z_{t_{n}}\cdot \mathbf{n})(X^r_{t_n,k})\vert^2\right],
	\end{aligned}
\end{equation}
where $X^l_{t_n,k}$ denotes the $k$-th point among the $K_l$ points selected on $\partial\Omega_{l}$ at time $t_n$ and $X^r_{t_n,k}$ denotes the $k$-th point among the $K_r$ points selected on $\partial\Omega_{r}$ at time $t_n$.  The parameters are chosen as $T = 1.0$, $f = 0$, $ \nu = 0.025 $, $ N = 50 $, $K_l = K_r = 100$, $\Delta t = 0.02$, $\alpha_1=\alpha_2=\alpha_3=0.01$ and $u_{in}=3$.
  Similarly, We choose to use the time adaptive approach \RNum{2} mentioned in \cite{wight2020solving} to improve accuracy and save training time. At $T = 1.0$, the streamline is shown in Figure \ref{fig:fig14} with $\nu = 0.025$, which is consistent with the result obtained by traditional numerical methods.
 
\subsection{Cahn-Hilliard equation}
\label{sec06:sec02}
\begin{figure}[htp]
	\centering
	\subfigure[Relative $L^2$ errors for 2D]{		 
		\includegraphics[scale=0.3]{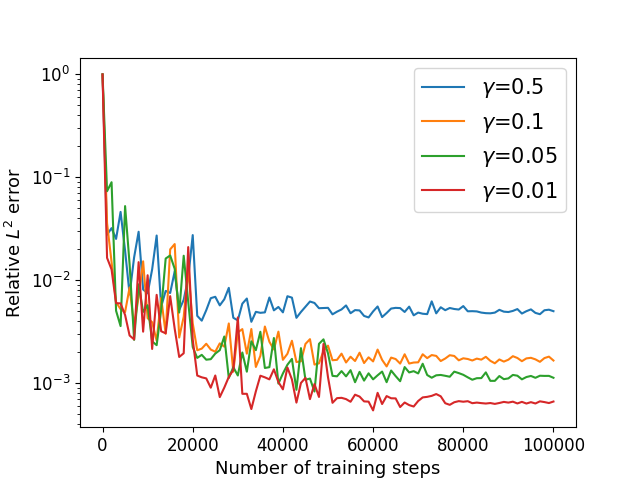}}
	\subfigure[Relative $L^2$ errors for 50D]{
		\includegraphics[scale=0.3]{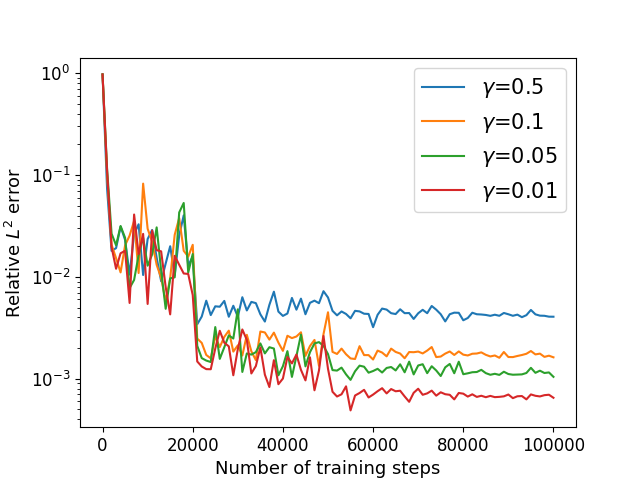}}
	\subfigure[Relative $L^2$ errors for 100D]{
		\includegraphics[scale=0.3]{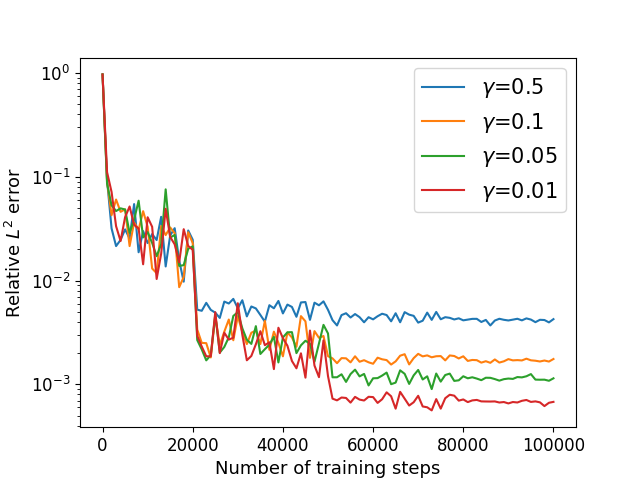}}
	\subfigure[Training losses for 2D]{		 
		\includegraphics[scale=0.3]{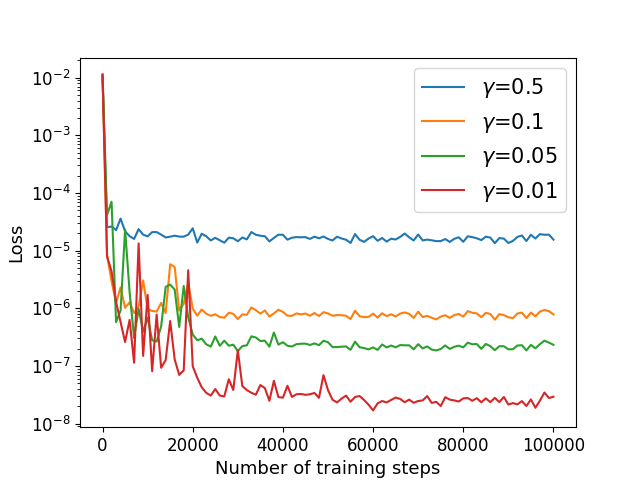}}
	\subfigure[Training losses for 50D]{
		\includegraphics[scale=0.3]{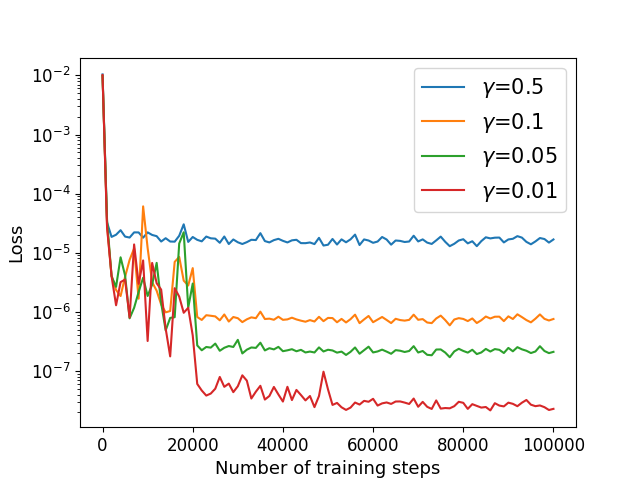}}
	\subfigure[Training losses for 100D]{
		\includegraphics[scale=0.3]{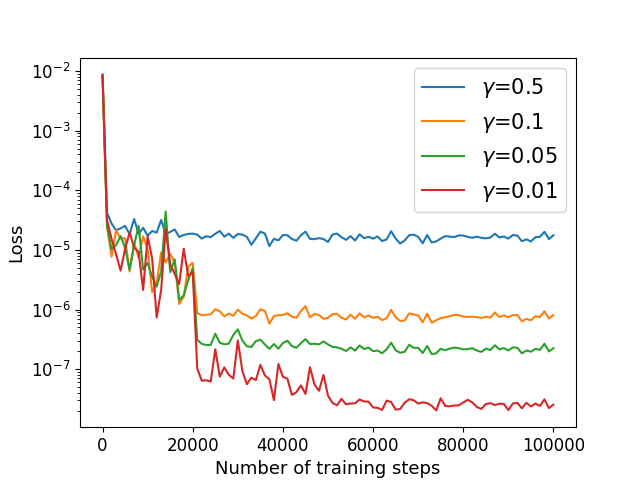}}
	\caption{Relative $L^{2}$ errors and training losses of the Algorithm \ref{alg:Algorithm 2} for \eqref{eq:ch2de} with different  parameter $\gamma$ in different dimension.}
	\label{fig:fig7}
\end{figure}

In this section, we consider the Cahn-Hilliard equation \eqref{eq:CH} 
in $(t,x)=(0,T]\times[-1,1]^d$ with initial condition
\begin{equation}\label{eq:ch2deinit}
	\begin{aligned} 
		& \phi(0,x) = \cos\left(\frac{\pi}{\sqrt{d}}\sum_{i=1}^{d}x_i\right).\\
	\end{aligned}
 \end{equation}
The exact solution is given by
\begin{equation}\label{eq:ch2de}
	\begin{aligned} 
		& \phi(t,x) = e^{-t}\cos\left(\frac{\pi}{\sqrt{d}}\sum_{i=1}^{d}x_i\right).\\
	\end{aligned}
\end{equation}
The parameters are taken as $L_{d}=$ 5E-4, $T=0.1$, $N=10$, $\delta=\Delta t=0.01$ and  $\alpha_1 = 0.01$. We estimate $\phi$ using Algorithm \ref{alg:Algorithm 2} with different parameter $\gamma$ in different dimension. The numerical results of the  errors for $\phi$ with different $\gamma$ and $S$ are recorded in Table \ref{Table:5}. Training processes in different dimension are shown in Figure \ref{fig:fig7}. For a fixed dimension, when $\gamma$ decreases, the relative $L^2$ error and training losses decrease. Our method is not sensitive to parameters $\gamma$ and $S$, and the training time of our method increases linearly with the dimension $d$, while the accuracy does not decrease. It works for the problem with high-order derivatives in high dimensions, which does not make the training difficult.
\begin{table}[htp]
	\begin{center} 
		\begin{tabular}{|l|l|l|l|l|l|l|l|l|}
				\hline
			d  &$\gamma$& S & $\vert\vert e\vert\vert_{L^{\infty}} $ &   $\vert\vert e\vert\vert_{L^2}$ & time      \\
			\hline 
			\multirow{4}{*}{2} & 0.5& 0.5 & 3.25E-2 & 5.02E-3 & \multirow{4}{*}{1200s} \\
			\cline{2-5} 
			& 0.1& 0.1 & 2.56E-3 & 1.67E-3 &\\
			\cline{2-5}
			& 0.05&  0.05 & 2.08E-3 &  1.13E-3  &\\
			\cline{2-5}
			& 0.01&  0.01 & 1.99E-3 & 6.67E-4  &\\
			\hline
			\multirow{4}{*}{50} & 0.5& 0.5 & 1.08E-2 & 4.06E-3 & \multirow{4}{*}{3200s} \\
			\cline{2-5} 
			& 0.1& 0.1 & 3.78E-3 & 1.65E-3   &\\
			\cline{2-5}
			& 0.05&  0.05 & 3.48E-3 &  1.04E-3 &\\
			\cline{2-5}
			& 0.01&  0.01 & 3.63E-3 & 6.49E-4 &\\
			\hline
			\multirow{4}{*}{100} & 0.5& 0.5 & 1.62E-2 & 4.25E-3 & \multirow{4}{*}{5200s} \\
			\cline{2-5}  
			& 0.1 & 0.1 & 6.46E-3 &  1.76E-3 & \\
			\cline{2-5}
			& 0.05 &  0.05 & 2.38E-3  & 1.15E-3 &\\
			\cline{2-5}
			& 0.01 &  0.01 & 4.06E-3 &  6.79E-4 &\\
			\hline
		\end{tabular}
		\caption{Relative $L^{\infty}$ and  $L^2$  errors  for \eqref{eq:ch2de} performed by the Algorithm \ref{alg:Algorithm 2} with $\gamma=0.5,\ 0.1,\ 0.05,\ 0.01$ for dimension $d=2,\ 50,\ 100$. }\label{Table:5}
	\end{center}
\end{table}

\begin{figure}[htp]
	\centering
	\subfigure[Relative $L^2$ errors for 2D]{		 
		\includegraphics[scale=0.3]{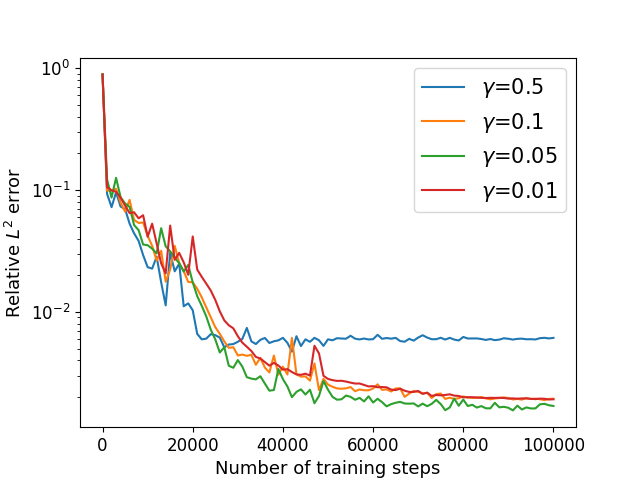}}
	\subfigure[Relative $L^2$ errors for 50D]{
		\includegraphics[scale=0.3]{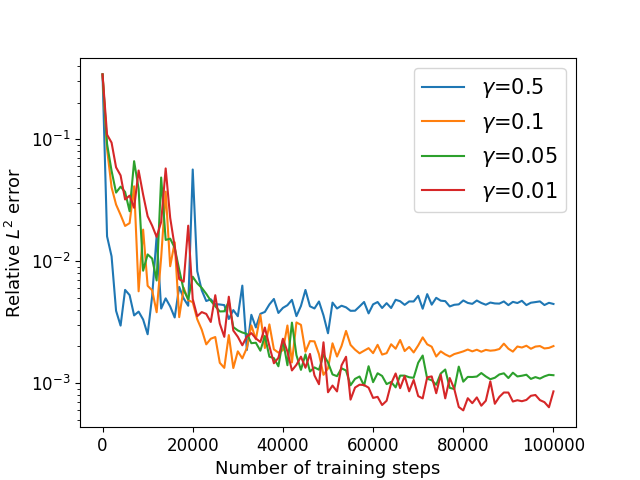}}
	\subfigure[Relative $L^2$ errors for 100D]{
		\includegraphics[scale=0.3]{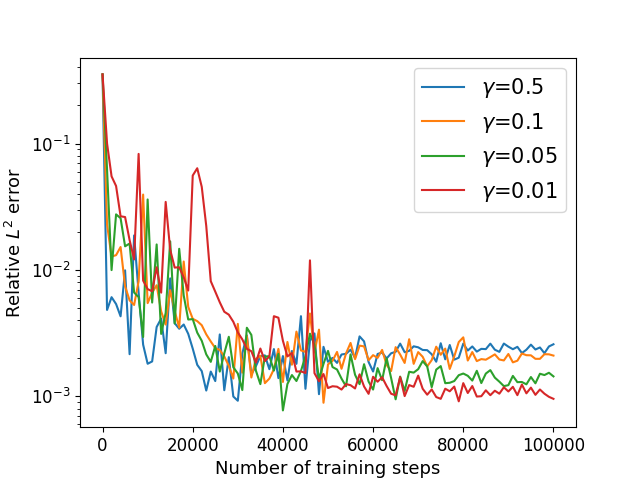}}
	\subfigure[Training loss for 2D]{		 
		\includegraphics[scale=0.3]{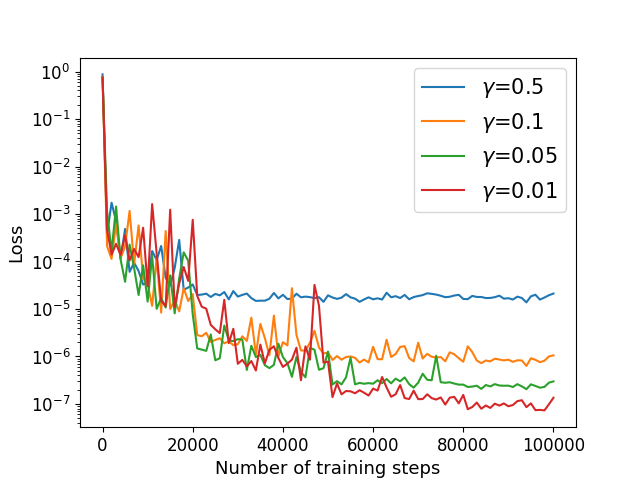}}
	\subfigure[Training loss for 50D]{
		\includegraphics[scale=0.3]{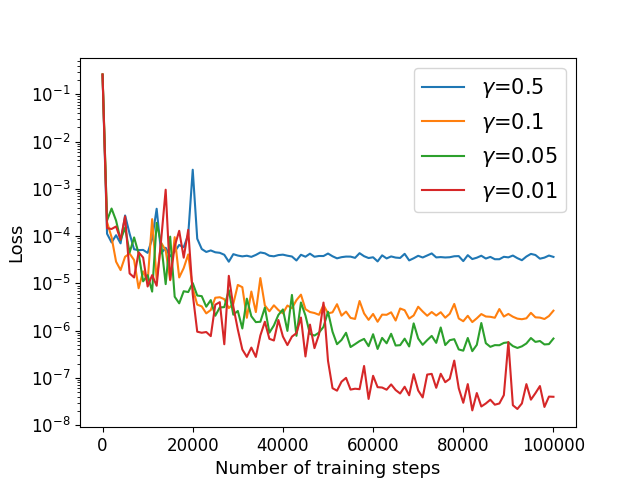}}
	\subfigure[Training loss for 100D]{
		\includegraphics[scale=0.3]{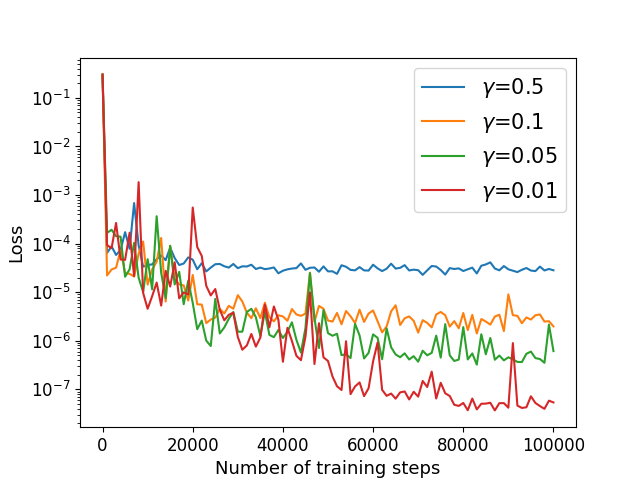}}
	\caption{Relative $L^{2}$ errors and training losses of the algorithm in Section \ref{sec04:sec03} for \eqref{eq:ch2de} with different parameter $\gamma$ in different dimension.} 
	\label{fig:fig8} 
\end{figure}

We consider the Cahn-Hilliard equation \eqref{eq:CH} with exact solution \eqref{eq:ch2de} defined in $\Omega=\{x:\vert x\vert < 1\} $ satisfying the mixed boundary conditions \eqref{eq:dnboundary1} on $\partial\Omega=\{x:\vert x\vert =1\}$, where $h(t,x)$ and $q(t,x)$ are given by the exact solution.
In this case, we choose $\alpha_1 = 1$. We utilize the algorithm in Section \ref{sec04:sec03} to solve $\phi$ with different parameter $\gamma$ in different dimension. 
The numerical results of the  errors for $\phi$ with different $\gamma$ and $S$ are recorded in Table \ref{Table:6} and the training processes are shown in Figure \ref{fig:fig8}. Our method works for boundary value problem in high dimension. 
\begin{table}[htp]
	\begin{center} 
		\begin{tabular}{|l|l|l|l|l|l|l|l|l|}
				\hline
			d  &$\gamma$& S & $\vert\vert e\vert\vert_{L^{\infty}} $ &   $\vert\vert e\vert\vert_{L^2}$ & time        \\
			\hline 
			\multirow{4}{*}{2} & 0.5& 0.5 & 1.69E-2 & 6.12E-3 & \multirow{4}{*}{1800s} \\
			\cline{2-5} 
			& 0.1& 0.1 & 4.52E-3 & 1.94E-3 & \\
			\cline{2-5}
			& 0.05&  0.05 & 3.64E-3 &  1.69E-3&  \\
			\cline{2-5}
			& 0.01&  0.01 & 3.72E-3 & 1.92E-3 &  \\
			\hline
			\multirow{4}{*}{50} & 0.5& 0.5 & 2.40E-2 & 4.46E-3 & \multirow{4}{*}{3800s} \\
			\cline{2-5} 
			& 0.1& 0.1 & 1.13E-2 & 2.01E-3 &   \\
			\cline{2-5}
			& 0.05&  0.05 & 5.77E-3 &  1.16E-3 & \\
			\cline{2-5}
			& 0.01&  0.01 & 7.33E-3 & 8.55E-4 & \\
			\hline
			\multirow{4}{*}{100} & 0.5& 0.5 & 1.14E-2 & 2.58E-3 &\multirow{4}{*}{6000s}\\
			\cline{2-5}  
			& 0.1 & 0.1 & 9.52E-3 &  2.10E-3 & \\
			\cline{2-5}
			& 0.05 &  0.05 & 1.10E-2  & 1.44E-3 & \\
			\cline{2-5}
			& 0.01 &  0.01 & 4.98E-2 &  9.56E-4 &\\
			\hline
		\end{tabular}
		\caption{Relative $L^{\infty}$ and  $L^2$ errors  for \eqref{eq:ch2de} with the mixed boundary conditions performed by the algorithm in Section \ref{sec04:sec03}, with $\gamma=0.5,\ 0.1,\ 0.05,\ 0.01$ for dimension $d=2,\ 50,\ 100$.} \label{Table:6}
  
	\end{center}
\end{table}

\begin{figure}[htp]
	\centering 
	\subfigure[Relative $L^2$ errors]{		 
		\includegraphics[scale=0.45]{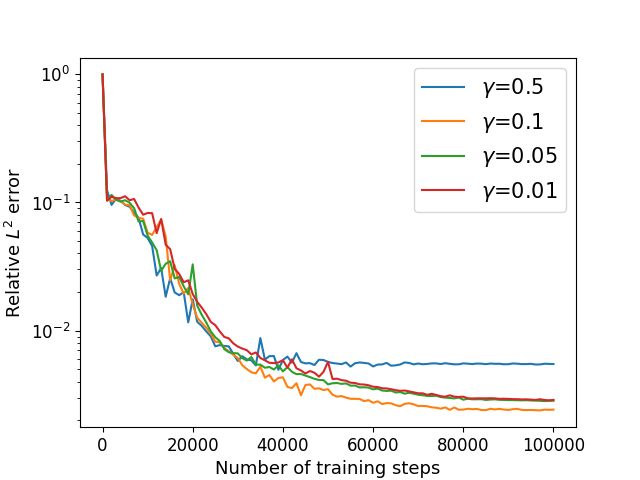}}
	\subfigure[Training loss]{		 
		\includegraphics[scale=0.45]{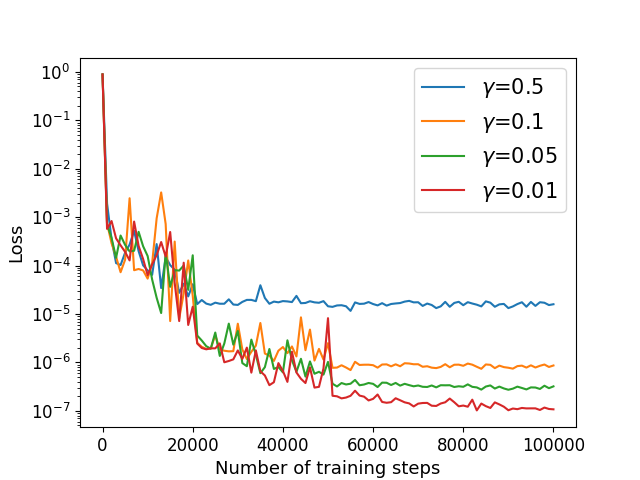}}
	\caption{Relative $L^2$ errors and training losses of the algorithm in Sections \ref{sec04:sec03} and \ref{sec04:sec04} for \eqref{eq:ch2de} with different parameter $\gamma$.}
	\label{fig:fig9}
\end{figure}

Next, we consider the Cahn-Hilliard equation \eqref{eq:CH} with exact solution \eqref{eq:ch2de} defined in $\Omega=(-\sqrt{2},\sqrt{2})^2$ with the periodic boundary condition \eqref{eq:pboundary1} and mixed boundary condition \eqref{eq:dnboundary1} on $\partial\Omega$, where the periods $I_i = 2\sqrt{2}, i=1,2$, where $h(t,x)$ and $q(t,x)$ are given by the exact solution. We choose $J=1$ and  other  parameters remain the same as previous example. We utilize the algorithm in Sections \ref{sec04:sec03} and \ref{sec04:sec04} to solve $\phi$ with different parameter $\gamma$. The numerical results of the errors for $\phi$ are recorded in Table \ref{Table:7} and Figure \ref{fig:fig9}. 
\begin{table}[htp]
	\begin{center} 
		\begin{tabular}{|l|l|l|l|l|l|l|l|l|}
            \hline
			$\gamma$& S & $\vert\vert e\vert\vert_{L^{\infty}}$ &  $\vert\vert e\vert\vert_{L^2}$ & time     \\
			\hline 
			0.5 & 0.5 & 9.98E-3 & 5.48E-3 & \multirow{4}{*}{3400s} \\
			\cline{1-4}
			0.1 & 0.1 & 4.86E-3 & 2.42E-3 &\\
			\cline{1-4}
			0.05& 0.05 & 7.86E-3 &  2.84E-3 & \\
			\cline{1-4}
			0.01& 0.01 & 7.39E-3 &  2.88E-3 &\\
			\hline
		\end{tabular}
		\caption{Relative $L^{\infty}$ and  $L^2$ errors for \eqref{eq:ch2de} with the mixed and periodic boundary condition performed by the algorithm in Sections \ref{sec04:sec03} and \ref{sec04:sec04}  with $\gamma=0.5,\ 0.1,\ 0.05,\ 0.01$. }\label{Table:7}
	\end{center}
\end{table}

\subsection{Cahn-Hilliard-Navier-Stokes equation}
\label{sec06:sec03}
\begin{figure}[htp]
	\centering 
	
	\subfigure[Relative $L^2$ errors]{		 
		\includegraphics[scale=0.45]{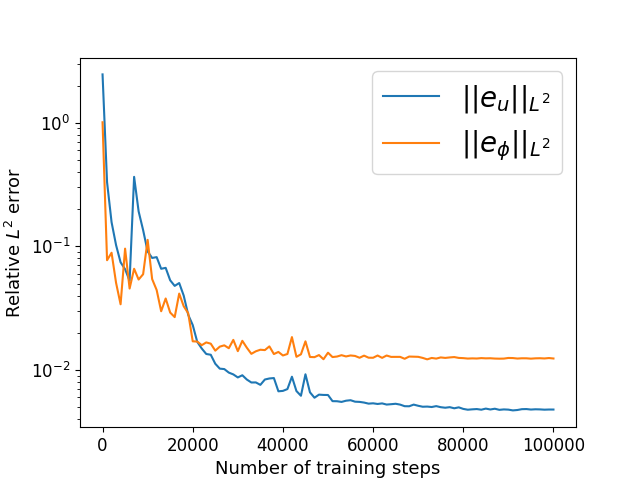}}
	\subfigure[Training loss]{		 
		\includegraphics[scale=0.45]{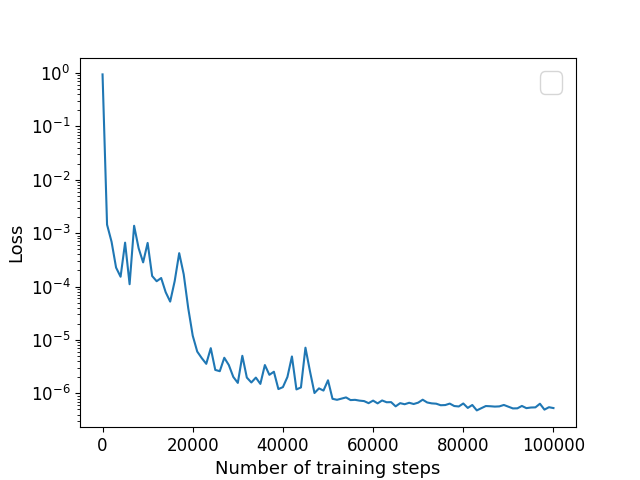}} 
	\caption{ Relative $L^2$ errors and training losses of the Algorithm \ref{alg:Algorithm 3} for \eqref{eq:CHNSexact}.}
	\label{fig:fig10}
\end{figure}

 In this section, we consider the coupled CHNS system  
\begin{equation}\label{eq:CHNSexact}
	\left\{
	\begin{aligned}
		& u_1(t,x)= -\cos(x_1)\sin(x_2)e^{-t},\\
		& u_2(t,x)= \sin(x_1)\cos(x_2)e^{-t},\\
		& p(t,x)=-\frac{1}{4}(cos(2x_1)+\cos(2x_2))e^{-2t} + c,\\
		& \phi(t,x)=\sin(x_1)\sin(x_2)e^{-t},
	\end{aligned}
	\right.
\end{equation}
for $(t,x) = (0,T]\times[0,2\pi]^2$ with the constant $c$ and initial condition
\begin{equation}\label{eq:CHNSexactinit}
	\left\{
	\begin{aligned}
		& u_1(0,x)= -\cos(x_1)\sin(x_2),\\
		& u_2(0,x)= \sin(x_1)\cos(x_2),\\
		& \phi(0,x)=\sin(x_1)\sin(x_2).\\
	\end{aligned}
	\right.
\end{equation} 
The parameters  are taken as $T = 0.1$,  $ N=5 $, $\delta=\Delta t =0.02$, $\nu = 1E-3$, $C = 1$, $L_d = 5E-4$, $\gamma = 0.01$, $S=0.0032$, $\alpha_1 = \alpha_2=\alpha_3 = 0.01$. The numerical results of the  errors for $\mathbf{u}(T,x)=(u_{1},u_2)^{T}$, $\phi$ and $\nabla p$ are recorded in Table \ref{Table:8} and Figure \ref{fig:fig10} shows the training process, where the Algorithm \ref{alg:Algorithm 3} in Section \ref{sec05} is implemented. It is easy to see that our method works for the coupled system.
\begin{table}[htp]
	\begin{center} 
		\begin{tabular}{|l|l|l|l|l|l|l|l|l|}
			\hline
			$\vert\vert e_{u_1}\vert\vert_{L^{\infty}}$ & $\vert\vert e_{u_2}\vert\vert_{L^{\infty}}$ & $\vert\vert e_{\phi}\vert\vert_{L^{\infty}}$ &$\vert\vert e_{u_1}\vert\vert_{L^{2}}$  &  $\vert\vert e_{u_2}\vert\vert_{L^{2}}$ &  $\vert\vert e_{\phi}\vert\vert_{L^{2}}$ & $\vert\vert e_{\nabla p}\vert\vert_{L^2} $& time    \\
			\hline 
		 	2.41E-2 & 2.30E-2 & 1.25E-2  &	1.31E-2 & 1.16E-2 & 4.78E-3 & 2.07E-1 & 1800s\\
			\hline 
		\end{tabular}
		\caption{Relative $L^{\infty}$ and  $L^2$ errors  for \eqref{eq:CHNSexact}  performed by the Algorithm \ref{alg:Algorithm 3}. }\label{Table:8}
	\end{center}
\end{table}

\begin{figure}[htp]
	\centering 
	\centering 
	\subfigure[t=0]{		 
		\includegraphics[scale=0.20]{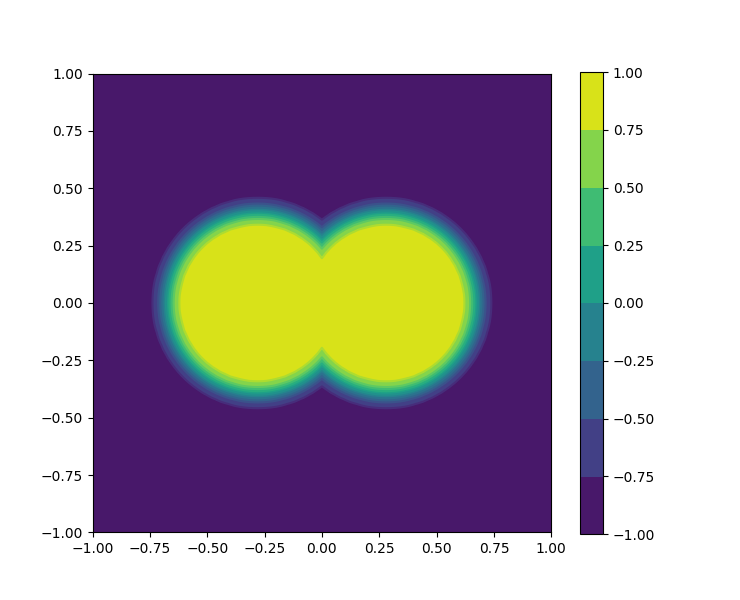}}
	\subfigure[t=0.2]{		 
		\includegraphics[scale=0.20]{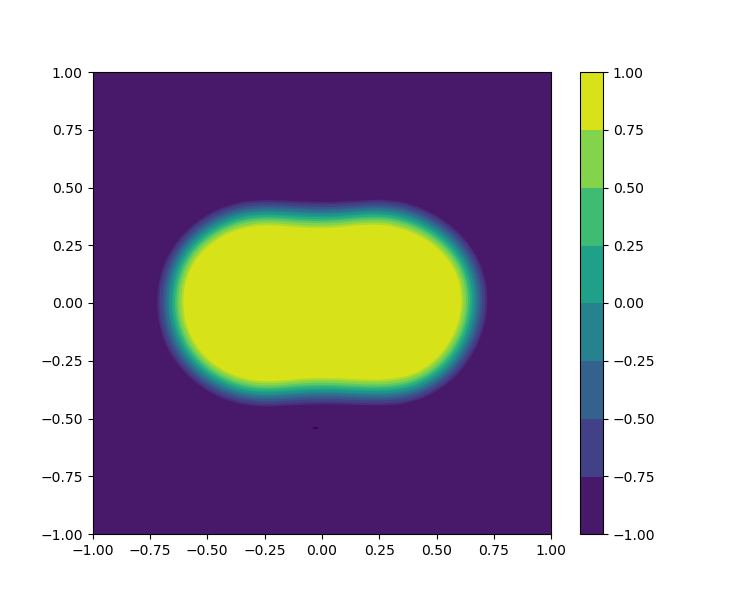}} 
  \subfigure[t=0.5]{		 
		\includegraphics[scale=0.20]{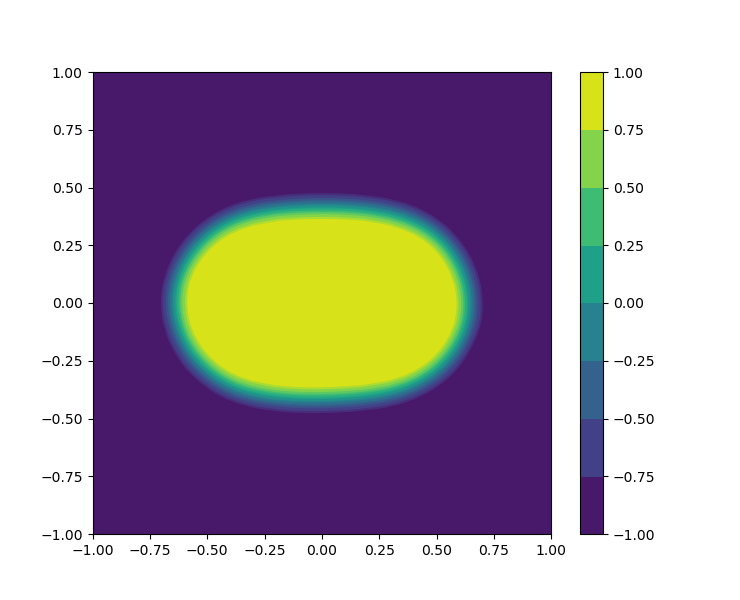}} 
  \subfigure[t=1.0]{		 
		\includegraphics[scale=0.20]{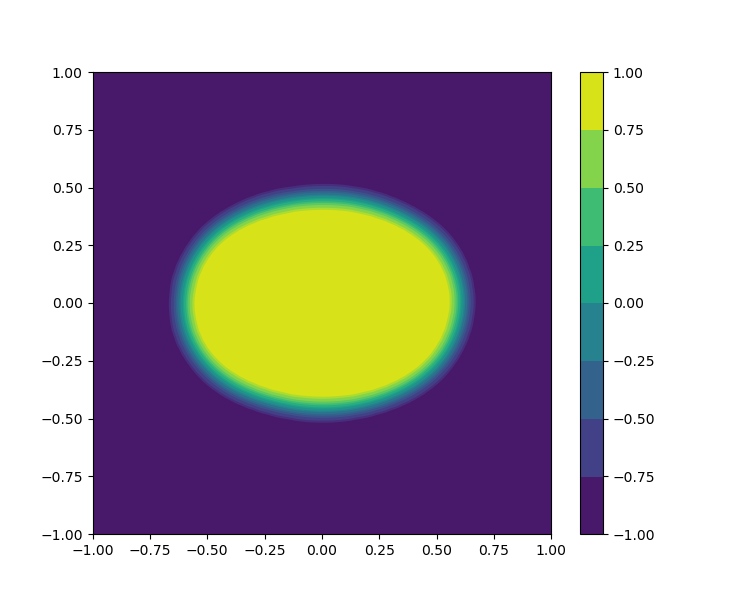}} 
  \subfigure[t=1.5]{		 
		\includegraphics[scale=0.20]{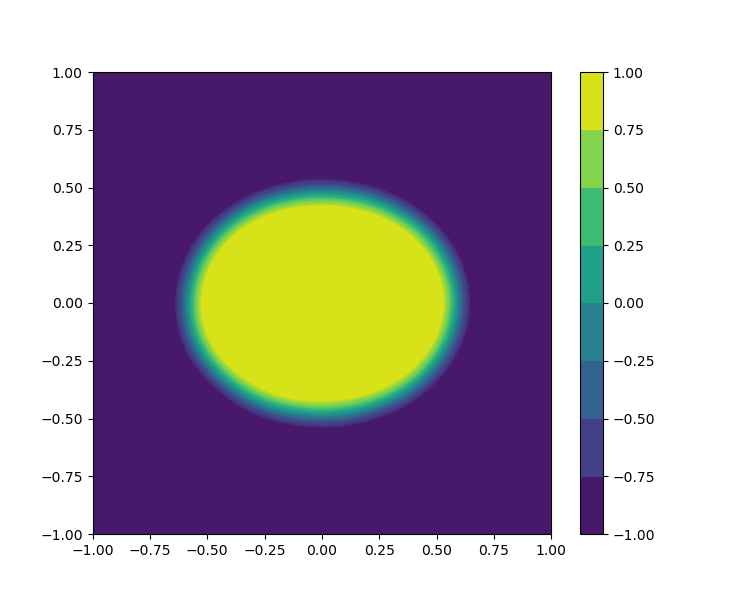}}
	\subfigure[t=2.0]{		 
		\includegraphics[scale=0.20]{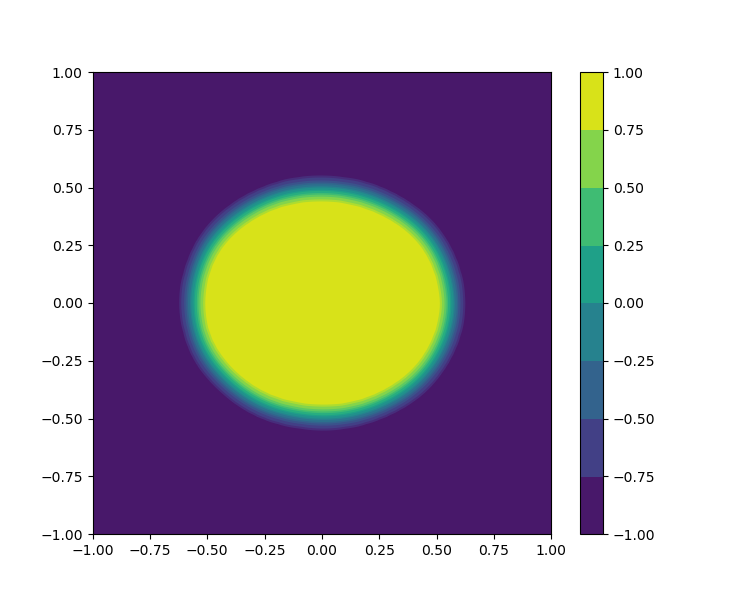}} 
  \subfigure[t=2.5]{		 
		\includegraphics[scale=0.20]{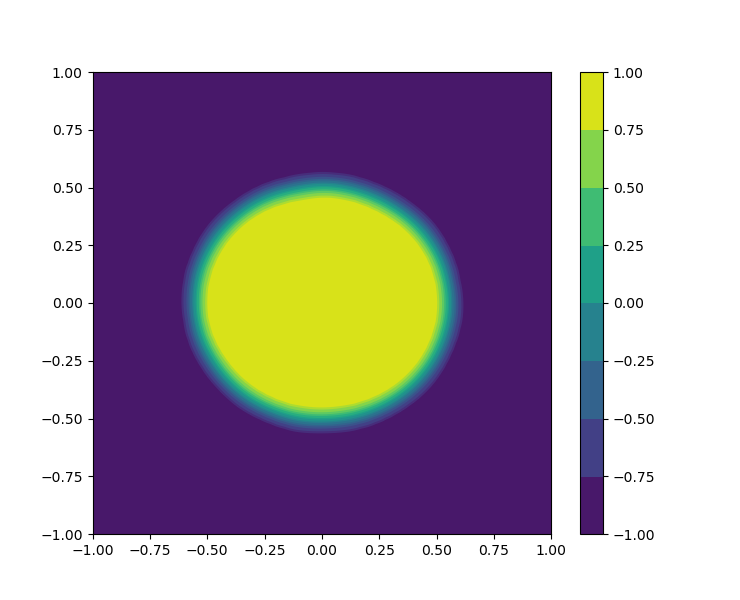}} 
  \subfigure[t=3.0]{		 
		\includegraphics[scale=0.20]{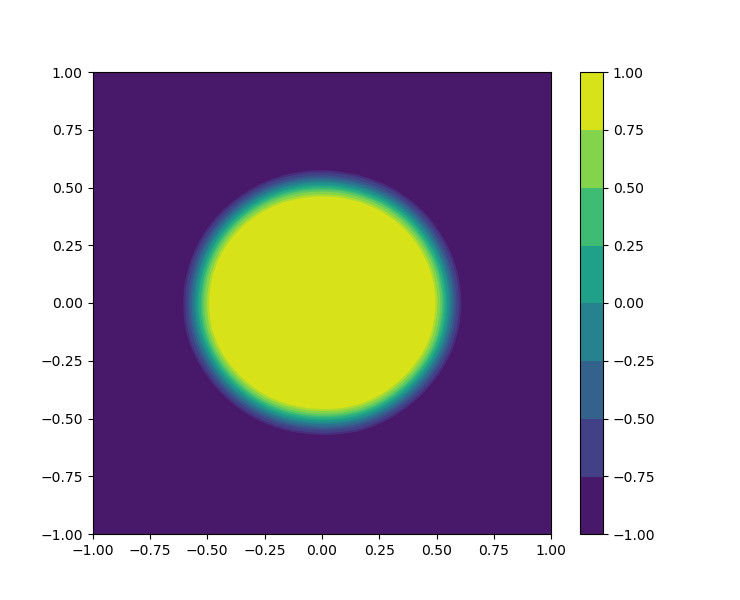}} 
	\caption{ Phase evolution at t = 0, 0.2, 0.5, 1.0, 1.5, 2.0, 2.5, 3.0.}
	\label{fig:fig15}
\end{figure}

Finally,  we study the interface problem modeled by the CHNS system. 
In this example,
we choose  the square domain $\Omega=[-1,1]^2$ and the parameters $L_d = 1$, $\nu = 1$, $C = 10$, $T=3$, $N = 300$, $\delta=\Delta t=0.01$, $\gamma=0.03$, $S=3.3$, $\alpha_1=\alpha_2=\alpha_3 = 0.01$. The 
initial conditions for $\phi$ and $\mathbf{u} = (u_1,u_2)^{T}$ is given
\begin{equation}\label{eq:interfacechns}
\left\{
    \begin{aligned}
    & \phi(0,x_1,x_2) = \max\left(\tanh\frac{r-R_1}{2\gamma},\tanh\frac{r-R_2}{2\gamma}\right),\\
    & \mathbf{u}(0,x_1,x_2) = 0,
    \end{aligned}
    \right.
\end{equation} 
where $r=0.4$, $R_1=\sqrt{(x_1-0.7r)^2+x_2^2}$, and $R_2=\sqrt{(x_1+0.7r)^2+x_2^2}$.
 The time adaptive approach \RNum{2} in \cite{wight2020solving} is employed to reduce the training time. According to the conservation of mass, we add the following loss term to the final loss function defined in Algorithm \ref{alg:Algorithm 3}
\begin{equation}\label{eq:mass}
	\begin{aligned} 
		& \frac{1}{N+1}\sum_{n=0}^{N}\left\vert\int_{\Omega}Y_{t_{n}}^\phi(X)dX-\int_{\Omega}g(X)dX\right\vert^2.
	\end{aligned}
\end{equation}
 The cosine and tanh functions are chosen as activation functions for the FNNs $\mathcal{U}_{\theta_1}$ and  $\mathcal{U}_{\theta_2}$, respectively.
 The evolution of the bubbles merging is visually shown in Figure \ref{fig:fig15}, which is coincide with the result in the literature.

\section{Conclusions and remarks}
\label{sec07}
In this article, we have presented the methods to obtain the numerical solutions of the incompressible Navier-Stokes equation, the Cahn-Hilliard equation and the CHNS system with different boundary conditions based on the Forward-Backward Stochastic Neural Networks. In particular, we utilize the modified Cahn-Hilliard equation that is derived from a widely used stabilized scheme for original Cahn-Hilliard, which can be diagonalized into a parabolic system. The FBSNNs are applied to this system, which works well for high-dimensional problem.  We demonstrate the performance of our algorithms on a variety of numerical experiments. In all numerical results, our methods are shown to be both stable and accurate. In the future work, we will study on how to make the training more efficiently and provide the theoretical analysis for our methods with some assumptions.

\section*{Declarations}
\begin{itemize}

\item \textbf{--Ethical Approval}\\
Not Applicable
 \item \textbf{--Availability of supporting data}\\
 The datasets generated during and/or analysed during the current study are available from the corresponding author on reasonable request.
 \item \textbf{--Competing interests}\\
 The authors have no relevant financial or non-financial interests to disclose. 
\item \textbf{--Funding}\\
This research is partially supported by the National key R \& D Program of China (No.2022YFE03040002) and the National Natural Science Foundation of China ( No.12371434).
\item \textbf{--Authors' contributions}\\
All authors contribute to the study conception and design. Numerical simulations are performed by Deng Yangtao. All authors read and approve the final manuscript.
\item\textbf{--Acknowledgments}\\
This research is partially supported by the National key R \& D Program of China (No.2022YFE03040002) and the National Natural Science Foundation of China ( No.12371434).
 \end{itemize}


\bibliographystyle{elsarticle-num}
\bibliography{journal}







\end{document}